\documentclass[12pt, a4paper]{amsart}
\usepackage{amscd, amssymb, pb-diagram}
\usepackage{amsthm}

\def\End{\operatorname{End}}
\def\newspan{\operatorname{span}}
\def\range{\operatorname{range}}

\def\ker{\operatorname{ker}}

\def\id{\operatorname{id}}
\def\iso{\operatorname{iso}}
\def\C{\mathbb{C}}
\def\R{\mathbb{R}}
\def\N{\mathbb{N}}

\newtheorem{thm}{Theorem}[section]
\newtheorem{cor}[thm]{Corollary}
\newtheorem{lemma}[thm]{Lemma}
\newtheorem{prop}[thm]{Proposition}

\theoremstyle{definition}
\newtheorem{definition}[thm]{Definition}

\theoremstyle{remark}

\newtheorem{example}[thm]{Example}

\numberwithin{equation}{section}

\begin{document}

\title[Partial-isometric crossed products]{Partial-isometric crossed products\\
by semigroups of endomorphisms}

\author{Janny Lindiarni}
\author{Iain Raeburn}
\address{School of Mathematical and Physical Sciences\\
University of Newcastle\\
NSW 2308\\
Australia} \email{janny, iain@frey.newcastle.edu.au}
\date{16 May 2002}
\thanks{This research was supported by the Australian Research Council.}

\begin{abstract}
Let $\Gamma^+$ be the positive cone in a totally ordered abelian
group $\Gamma$, and let $\alpha$ be an action of $\Gamma^+$ by
endomorphisms of a $C^*$-algebra $A$. We consider a new kind of
crossed-product $C^*$-algebra $A\times_\alpha\Gamma^+$, which is
generated by a faithful copy of $A$ and a representation of
$\Gamma^+$ as partial isometries. We claim that these crossed
products provide a rich and tractable family of Toeplitz algebras
for product systems of Hilbert bimodules, as recently studied by
Fowler, and we illustrate this by proving detailed structure
theorems for actions by forward and backward shifts.
\end{abstract}

\maketitle

\section{Introduction}

Let $\Gamma$ be a totally ordered abelian group with positive cone
$\Gamma^+$, and consider an action $\alpha : \Gamma^+ \rightarrow
\End A$ of $\Gamma^+$ by endomorphisms of a $C^*$-algebra $A$. We
study covariant representations $(\pi , V)$ of the system $(A,
\Gamma^+, \alpha)$ in which the endomorphisms $\alpha_s$ are
implemented by partial isometries $V_s$, and the corresponding
crossed-product $C^*$-algebra $A\times_\alpha \Gamma^+$ which is
generated by a universal covariant representation. We think these
\emph{partial-isometric crossed products} are likely to be of
interest for several reasons.

Our first motivation comes from the analogous covariant
representation theory in which the elements of $\Gamma^+$ are
implemented by isometries. To avoid confusion, we shall refer to
these as \emph{covariant isometric representations} and the
corresponding crossed products $A \times_{\alpha}^{\iso} \Gamma^+$
as \emph{isometric crossed products}. Isometric crossed products
by the semigroup $\mathbb{N} = \mathbb{Z}^+$ were first used to
give a model for the Cuntz algebra $\mathcal{O}_n$ \cite{cuntz,
paschke, stacey, bkr}. Subsequently, various authors considered
the action $\tau$ of $\Gamma^+$ by right translation on a
distinguished subalgebra $B_{\Gamma^+}$ of $\ell^{\infty}
(\Gamma^+)$, and used $B_{\Gamma^+} \times_{\tau}^{\iso} \Gamma^+$
as a model for the Toeplitz algebra $\mathcal{T}(\Gamma)$ . This
program has also been effective for the more general quasi-lattice
ordered groups, such as $\N^k$ \cite{lr}. More recently, isometric
crossed products by actions of $\mathbb{N}^k$ have proved to be
useful models for Hecke algebras arising in number theory (see
\cite{lr2, laca, larsen}, for example).

The theory of isometric crossed products, however, yields no
information about some systems $(A, \Gamma^+, \alpha)$. For
example, consider the action $\sigma$ of $\mathbb{N}$ by left
translations on
\[\mathbf{c}_0 := \{f: \mathbb{N} \to \mathbb{C}: f(n)
\rightarrow 0 \mbox{ as } n \rightarrow \infty\}.\] Every
covariant isometric representation $(\pi , V)$ of $(\mathbf{c}_0,
\mathbb{N}, \sigma)$ satisfies \[\pi (f) = (V^*)^n \pi (\sigma_n
(f)) V^n,\] and since $\sigma_n (f) \rightarrow 0$ as $n
\rightarrow \infty$ for every $f \in \mathbf{c}_0$, this is only
possible if $\pi$ is identically zero. Thus $\mathbf{c}_0
\times_{\sigma}^{\iso} \mathbb{N} = \{0\}$. Every system $(A,
\Gamma^+, \alpha)$, on the other hand, admits covariant
partial-isometric representations $(\pi , V)$ in which $\pi$ is
faithful, and hence the partial-isometric crossed product contains
full information about the system.

Our second motivation concerns the Toeplitz algebras of Hilbert
bimodules \cite{fr}. In Pimsner's original investigations of
Hilbert bimodules \cite{pimsner}, a key example was provided by an
endomorphism $\alpha$ of a $C^*$-algebra $A$, and the
Cuntz-Pimsner algebra of this bimodule is the isometric crossed
product $A \times_{\alpha}^{\iso} \mathbb{N}$. The Toeplitz
algebra of this bimodule, on the other hand, is our
partial-isometric crossed product $A\times_\alpha\mathbb{N}$.
Fowler has recently considered product systems of Hilbert
bimodules over more general semigroups, and studied the Toeplitz
algebras of these product systems. In particular he has identified
conditions under which the results of \cite{fr} carry over to his
new family of Toeplitz algebras \cite{f}. Important examples of
product systems are provided by endomorphic actions of semigroups,
and our partial-isometric crossed products are Toeplitz algebras
to which Fowler's results apply. Because of their concrete nature,
partial-isometric crossed products form a particularly tractable
family of Toeplitz algebras, and even for the systems
$(B_{\mathbb{N}}, \mathbb{N}, \tau)$ and $(B_{\mathbb{N}},
\mathbb{N}, \sigma)$ the partial-isometric crossed products have
rich structure. Thus our results confirm that there is a lot of
interesting information lying between a Cuntz-Pimsner algebra and
its Toeplitz algebra extension. In particular, for the Hilbert
bimodules associated to the endomorphisms $\tau_1$ and $\sigma_1$
of $B_{\mathbb{N}}$, there are many distinct relative
Cuntz-Pimsner algebras as in \cite{ms, fmr}.

A third point of interest lies in the form of our structure
theorems for crossed products. Associated to any pair of ideals
$I$, $J$ in a $C^*$-algebra $B$ is a commutative diagram
\begin{equation*}
\begin{diagram}\dgARROWLENGTH=0.8\dgARROWLENGTH
\node{} \node{0} \arrow{s} \node{0} \arrow{s} \node{0} \arrow{s}
\\
\node{0} \arrow{e} \node{I \cap J} \arrow{e} \arrow{s} \node{J}
\arrow{e} \arrow{s} \node{J/(I \cap J)} \arrow{e} \arrow{s}
\node{0}
\\
\node{0} \arrow{e} \node{I} \arrow{e} \arrow{s} \node{B} \arrow{e}
\arrow{s} \node{B/I} \arrow{e} \arrow{s} \node{0}
\\
\node{0} \arrow{e} \node{I/(I \cap J)} \arrow{e} \arrow{s}
\node{B/J} \arrow{e} \arrow{s} \node{B/(I+J)} \arrow{e} \arrow{s}
\node{0}
\\
\node{} \node{0} \node{0} \node{0}
\end{diagram}
\end{equation*}
in which all the rows and columns are exact. From this, it follows
easily that $B$ has a composition series $0 \leq I \cap J \leq I+J
\leq B$ of ideals with subquotients
\[
I \cap J,\quad (I+J)/(I \cap J) \cong (I/ (I \cap J)) \oplus (J/
(I \cap J))\quad\text{and}\quad B/(I+J),
\]
and it is often helpful to have a structure theorem which
identifies these subquotients in familiar terms. Here, though, we
can do more: we can identify in familiar terms the four extensions
which make up the outside square. Thus, for example,
Theorem~\ref{diagramcNtau} describes $B = B_{\mathbb{N}}
\times_{\tau} \mathbb{N}$, and identifies both the right-hand and
bottom exact sequences as the extension of $C (\mathbb{T})$ by
$\mathcal{K}$ provided by the Toeplitz algebra $\mathcal{T}
(\mathbb{Z})$.

\smallskip

We begin with a preliminary section containing background material
about power partial isometries, Toeplitz algebras and isometric
crossed products, and Hilbert bimodules. In \S\ref{sectpireps}, we
discuss representations of totally ordered semigroups by partial
isometries, and analyse the $C^*$-algebras generated by two
semigroups of truncated shifts.

In \S\ref{sect2}, we discuss covariant partial-isometric
representations and the partial-isometric crossed product. We show
that every action $\alpha:\Gamma^+\to \End A$ admits covariant
partial-isometric representations in which $A$ acts faithfully
(Example~\ref{repwithAfaithful}), and show how the results of
Fowler \cite{f} allow us to identify the covariant
partial-isometric representations $(\pi,V)$ of
$(A,\Gamma^+,\alpha)$ for which the associated representation
$\pi\times V$ of $A\times_\alpha\Gamma^+$ is faithful
(Theorem~\ref{faithful}).

In \S\ref{sect3}, we give a structure theorem for the crossed
product of the system $(B_{\Gamma^+},\Gamma^+,\tau)$ arising in
the analysis of Toeplitz algebras. While many of our observations
work for arbitrary totally ordered abelian groups, the main
Theorem~\ref{Bgammatau} concerns subsemigroups of $\mathbb{R}^+$.
In Theorem~\ref{diagramcNtau}, we obtain more detailed information
for the semigroup $\N$. In the last section, we consider the
action $\sigma$ of $\N$ by right translation on $B_\N$. Although
the structure of $B_\N\times_\sigma \N$ is quite a bit simpler
than that of $B_\N\times_\tau\N$, it is still a good deal more
complicated than that of $B_\N\times_\sigma^{\iso} \N=\{0\}$ (see
Theorem~\ref{csigmaN}).

\subsection*{Conventions} Throughout this paper, $\Gamma$ will be a
totally ordered abelian group with positive cone $\Gamma^+$;
sometimes $\Gamma$ will be a subgroup of $\mathbb{R}$ but if so,
we shall say so. Our main examples are the additive semigroup $\N$
(which for us always contains $0$) and other subsemigroups of the
additive group $\R$; we therefore use additive notation in
$\Gamma$, so that the identity is $0$ and $\Gamma^+=\{s\in \Gamma:
x\geq 0\}$. A subgroup $I$ of $\Gamma$ is an \emph{order ideal} if
$x\in I$ and $0\leq y\leq x$ imply $y\in I$. We say that $\Gamma$
is \emph{simple} if it has no nontrivial order ideals; standard
theorems say that $\Gamma$ is simple if and only if it is
\emph{archimedean} in the sense that there exists $x\in \Gamma^+$
such that $\{y\in \Gamma^+:y\leq nx\text{ for some $n\in \N$}\}$
is all of $\Gamma^+$, and hence if and only if $\Gamma$ is order
isomorphic to a subgroup of $\R$ \cite{fuchs}.

We denote by $\mathcal{K}(H)$ the $C^*$-algebra of compact
operators on a Hilbert space $H$. We write $\bar\pi$ for the
extension of a non-degenerate representation $\pi:A\to B(H)$ to
the multiplier algebra $M(A)$; similarly, if $\alpha$ is an
extendible endomorphism of a $C^*$-algebra $A$, in the sense that
there is an approximate identity $\{a_i\}$ such that $\alpha(a_i)$
converges strictly to a projection in $M(A)$ \cite[\S2]{adji2}, we
write $\bar \alpha$ for the extension of $\alpha$ to an
endomorphism of $M(A)$.

\section{Preliminaries}\label{sect1}

\subsection{Power partial isometries}

A \emph{partial isometry} $V$ on a Hilbert space $H$ is an
operator on $H$ such that $\|Vh\|=\|h\|$ for $h \in (\ker
V)^{\perp}$. A bounded operator $V$ is a partial isometry if and
only if $VV^*V=V$, and then the operators $V^*V$ and $VV^*$ are
the orthogonal projections on the initial space $(\ker V)^{\perp}$
and range $VH$ respectively. An element $v$ of a $C^*$-algebra $A$
is called a partial isometry if $vv^*v=v$; it then becomes a
partial isometry in the usual sense whenever we represent $A$ on
Hilbert space.

The product of two partial isometries is not in general a partial
isometry: for example, if $\{e_1, e_2\}$ is an orthonormal set and
$P$, $Q$ are the orthogonal projections on $\newspan\{e_1\}$,
$\newspan\{e_1+e_2\}$, then $e_1$ is a unit vector in $(\ker
QP)^\perp$ but $QPe_1=(e_1+e_2)/2$ does not have norm $1$. Since
we are interested in semigroups of partial isometries, we need to
know when a product of partial isometries is a partial isometry.
The answer is well-known, and is proved, for example, in
\cite[Lemma~2]{hw}.

\begin{prop}\label{tool}
Let $S$ and $T$ be partial isometries. Then $ST$ is a partial
isometry if and only if $S^*S$ commutes with $TT^*$.
\end{prop}

A partial isometry $v$ is a \emph{power partial isometry} if $v^n$
is a partial isometry for $n \in \mathbb{N}$. Every isometry and
coisometry is a power partial isometry, and the other key examples
are the truncated shifts described in Example~\ref{ppisometryJk}
--- indeed, every power partial isometry is a direct sum of an
isometry, a coisometry and some truncated shifts \cite{hw}.

\begin{example}\label{ppisometryJk}
Let $\{e_i: 1 \leq i \leq k+1\}$ be the usual orthonormal basis
for $\mathbb{C}^{k+1}$, and consider the \emph{truncated shift}
$J_k = \sum_{j=1}^k e_{j+1} \otimes \overline{e}_j$ on
$\mathbb{C}^{k+1}$. Since the $e_i \otimes \overline{e}_j$ are
matrix units, we have $J_k^* = \sum_{j=2}^{k+1} e_{j-1} \otimes
\overline{e}_j$, and can check that $J_k$ is a power partial
isometry satisfying $J_k^{k+1} = 0$. The $C^*$-subalgebra $C^*
(J_k)$ of $B(\C^{k+1})$ generated by $J_k$ contains $e_i \otimes
\overline{e}_j = (J_k^*)^{j-1} J_k^k (J_k^*)^k J_k^{i-1}$, and
hence is all of $B (\mathbb{C}^{k+1})$.
\end{example}

The $C^*$-algebra $C^* (\bigoplus_{k=1}^\infty J_k)$ is universal
among $C^*$-algebras generated by a power partial isometry
\cite[Theorem 1.3]{hr}. We shall need to know about $C^*
(\bigoplus_{k=1}^n J_k)$.

\begin{lemma}\label{dirsumJk}
The $C^*$-algebra $C^*(\bigoplus_{k \leq n} J_k)$ is isomorphic to
$\bigoplus_{k \leq n} M_{k+1} (\mathbb{C})$.
\end{lemma}

\begin{proof}
We prove this by induction on $n$. For $n=1$, we have
$C^*(J_1)\cong B(\C^2)$ by Example~\ref{ppisometryJk}. Suppose
$n>1$ and we know the result for $n-1$. Then since $J_k^n=0$ for
$k<n$, we have $(\bigoplus_{k\leq n}J_k)^n=0\oplus
J_n^n=0\oplus(e_{n+1}\otimes \overline{e}_1)$. Since $C^*
(e_{n+1}\otimes \overline{e}_1) = B (\mathbb{C}^{n+1}) = C^*
(J_n)$, this implies that $(\bigoplus_{k< n}J_k)\oplus 0\in
C^*(\bigoplus_{k\leq n}J_k)$, and that
\[
\textstyle{C^*(\bigoplus_{k\leq n}J_k)=C^*(\bigoplus_{k<
n}J_k)\oplus B(\C^{n+1}),}
\]
which is isomorphic to $\bigoplus_{k\leq n} B(\C^{k+1})$ by the
inductive hypothesis.
\end{proof}

\subsection{Toeplitz algebras}\label{toepalg}

We let $\{\varepsilon_r: r \in \Gamma^+\}$ denote the usual
orthonormal basis of $\ell^2 (\Gamma^+)$. There is a
representation $T=T^\Gamma$ of $\Gamma^+$ by isometries on $\ell^2
(\Gamma^+)$ such that $T_s (\varepsilon_r) = \varepsilon_{r+s}$
for $r,s \in \Gamma^+$. The \emph{Toeplitz algebra} $\mathcal{T}
(\Gamma)$ is the $C^*$-subalgebra of $B (\ell^2 (\Gamma^+))$
generated by $\{T_s: s \in \Gamma^+\}$.

One way to analyse $\mathcal{T}(\Gamma)$ is by realising it as a
semigroup crossed product. If $\alpha$ is an action of $\Gamma^+$
by endomorphisms of a $C^*$-algebra $A$, then a \emph{covariant
isometric representation} of $(A,\Gamma^+,\alpha)$ consists of a
non-degenerate representation $\pi:A\to B(H)$ and an isometric
representation $V$of $\Gamma^+$ on $H$ such that
\[
\pi (\alpha_s (a)) = V_s \pi (a) V_s^*\ \text{ for $s \in
\Gamma^+$ and $a\in A$;}
\]
the semigroup crossed product of $A$ by $\alpha$ is by definition
universal for such representations. Here we call it the
\emph{isometric crossed product} and denote it by
$A\times_{\alpha}^{\iso} \Gamma^+$; we write $\pi\times^{\iso}V$
for the representation of $A\times_{\alpha}^{\iso} \Gamma^+$
associated to a covariant isometric representation $(\pi,V)$.

For $s\in \Gamma^+$, we denote by $1_s$  the characteristic
function of $\{t\in \Gamma^+:t\geq s\}$. Since $1_s 1_t = 1_{\max
\{s,t\}}$ and $1_s^* = 1_s$, $B_{\Gamma^+}:=\overline{\newspan}
\{1_t: t \in \Gamma^+\}$ is a $C^*$-subalgebra of $\ell^{\infty}
(\Gamma^+)$. The action $\tau$ of $\Gamma^+$ by right translation
on $\ell^\infty(\Gamma^+)$ satisfies $\tau_t(1_s)=1_{s+t}$, and
hence restricts to an action of $\Gamma^+$ by endomorphisms of
$B_{\Gamma^+}$. For every isometric representation $V$ of
$\Gamma^+$ on $H$, there is a representation $\pi_V$ of
$B_{\Gamma^+}$ such that $(\pi_V,V)$ is a covariant isometric
representation of $(B_{\Gamma^+},\Gamma^+,\tau)$; if each $V_s$ is
non-unitary, then $\pi_V\times^{\iso} V$ is a faithful
representation of $B_{\Gamma^+}\times_{\tau}^{\iso}\Gamma^+$
\cite[Theorem~2.4]{alnr}. In particular, the representation
$T=T^\Gamma:\Gamma^+\to B(\ell^2(\Gamma^+))$ induces an
isomorphism $\pi_T\times^{\iso}T$ of
$B_{\Gamma^+}\times_{\tau}^{\iso}\Gamma^+$ onto
$\mathcal{T}(\Gamma)$.

Since every unitary representation is in particular an isometric
representation, there is a canonical quotient map $\psi_T$ of
$\mathcal{T}(\Gamma)\cong
B_{\Gamma^+}\times_{\tau}^{\iso}\Gamma^+$ onto $C^*(\Gamma)\cong
C(\hat\Gamma)$. Murphy proved in \cite{murphy1} that the kernel of
$\psi_T$ is the commutator ideal $\mathcal{C}_\Gamma$ of
$\mathcal{T}(\Gamma)$. In \cite[Remark~3.3]{a}, this is deduced
from properties of isometric crossed products; in particular, it
is shown that $\mathcal{C}_\Gamma$  itself an isometric crossed
product $B_{\Gamma^+,\infty}\times_{\tau}^{\iso}\Gamma^+$. From
this description it follows easily that $\mathcal{C}_\Gamma$ is
generated by the elements $T_sT_s^*-T_tT_t^*$ for $s<t$.

Douglas proved that if $\Gamma$ is a subgroup of $\mathbb{R}$ then
$\mathcal{C}_\Gamma$ is a simple $C^*$-algebra, and Murphy proved
the converse (\cite[Theorem~4.3]{murphy1}; see also
\cite[page~1141]{alnr}). The following description of
$\mathcal{C}_\Gamma$ will be useful.

\begin{lemma}\label{comideal}
For every totally ordered abelian group $\Gamma$,
\begin{equation}\label{idcI}
\mathcal{C}_{\Gamma} = \overline{\newspan} \{T_r (1 - T_u T_u^*)
T_t^*: r,u,t \in \Gamma^+\}.
\end{equation}
\end{lemma}

\begin{proof}
Let $\mathcal{I}$ denote the right-hand side of (\ref{idcI}).
Since each $1-T_uT_u^*=[T_u^*,T_u]$ is a commutator,
$\mathcal{I}\subset \mathcal{C}_\Gamma$. On the other hand,
$\mathcal{C}_\Gamma$ is generated by the elements $T_s T_s^* - T_t
T_t^*=T_s (1 - T_{t-s} T_{t-s}^*) T_s^*$ for $s<t$, so it suffices
to prove that $\mathcal{I}$ is an ideal.

Since $\mathcal{T} (\Gamma)$ is generated by $\{T_s: s \in
\Gamma^+ \}$, we must show $T_s \mathcal{I} \subset \mathcal{I}$
and $T_s^* \mathcal{I} \subset \mathcal{I}$. We trivially have
$T_s \mathcal{I} \subset \mathcal{I}$. Let $r,u,t \in \Gamma^+$.
Then
\begin{align*}
T_s^* T_r (1 &- T_u T_u^*) T_t^* \\
&= \left\{ \begin{array}{ll}
T_{r-s} (1 - T_u T_u^*) T_t^* & \mbox{if $r \geq s$} \\
T_{s-r}^* (1 - T_u T_u^*) T_t^* & \mbox{if $r<s$}
\end{array}
\right. \\
&= \left\{ \begin{array}{ll}
T_{r-s} (1 - T_u T_u^*) T_t^* & \mbox{if $r \geq s$} \\
(T_{s-r}^* - T_{s-r}^* T_{s-r} T_{u-(s-r)} T_{u-(s-r)}^*
T_{s-r}^*) T_t^*
& \mbox{if $r<s$ and $s-r<u$} \\
0 & \mbox{otherwise}
\end{array}
\right. \\
&= \left\{ \begin{array}{ll}
T_{r-s} (1 - T_u T_u^*) T_t^* & \mbox{if $r \geq s$} \\
(1 - T_{u-(s-r)} T_{u-(s-r)}^*) T_{s-r+t}^* & \mbox{if $r<s$ and
$s-r<u$}
\\
0 & \mbox{otherwise}
\end{array}
\right.
\end{align*}
belongs to $\mathcal{I}$. Hence $T_s^* \mathcal{I} \subset
\mathcal{I}$, and $\mathcal{I}$ is an ideal, as required.
\end{proof}

\subsection{Hilbert bimodules}

A \emph{Hilbert bimodule} over a $C^*$-algebra $A$ is a right
Hilbert $A$-module $X$ together with a homomorphism $\phi : A
\rightarrow \mathcal{L} (X)$, which gives a left action of $A$ on
$X$: $a \cdot x := \phi (a) x$ for $a \in A$ and $x \in X$. A
\emph{Toeplitz representation} $(\psi , \pi)$ of $X$ in a
$C^*$-algebra $B$ consists of a linear map $\psi : X \rightarrow
B$ and a homomorphism $\pi : A \rightarrow B$ such that
\begin{align*}
\psi (x \cdot a) = \psi (x) \pi (a),\quad \psi (x)^* \psi (y) =
\pi (\langle x,y \rangle_A),\quad\text{and}\quad \psi (a \cdot x)
= \pi (a) \psi (x)
\end{align*}
for every $x \in X$ and $a \in A$.

The \emph{Toeplitz algebra} of $X$ is the $C^*$-algebra
$\mathcal{T}_X$ generated by the range of a universal Toeplitz
representation $(i_X, i_A)$ of $X$; then for every Toeplitz
representation $(\psi , \pi)$ of $X$ in a $C^*$-algebra $B$, there
is a homomorphism $\psi \times \pi$ of $\mathcal{T}_X$ into $B$
such that $(\psi \times \pi) \circ i_X = \psi$ and $(\psi \times
\pi) \circ i_A = \pi$. For every Hilbert bimodule $X$, there is
such a $C^*$-algebra and it is unique up to isomorphism
\cite[Proposition 1.3]{fr}.

\begin{definition}
Suppose $X = \{X_s : s \in \Gamma^+\}$ is a family of Hilbert
bimodules over a $C^*$-algebra $A$. We write $\phi_s: A
\rightarrow \mathcal{L} (X_s)$ for the left action of $A$ on
$X_s$. Fowler says that $X$ is a \emph{product system} over
$\Gamma^+$ if there is an associative multiplication on $X$
(strictly speaking, on the disjoint union of the $X_s$) such that
for $s,t \in \Gamma^+$, the map $(x,y) \mapsto xy$ extends to an
isomorphism of the Hilbert bimodule $X_s \otimes X_t$ onto
$X_{st}$ \cite[Definition 2.1]{f}. He also requires that $X_0 = A$
(with left and right actions defined by multiplication in $A$),
and that the multiplications in $X$ involving elements in $A =
X_0$ satisfy $ax = a \cdot x$ and $xa = x \cdot a$ for $a \in A$
and $x \in X_s$.
\end{definition}

\begin{definition}
A \emph{Toeplitz representation} $\psi$ of a product system $X$ in
a $C^*$-algebra $B$ is a map $\psi : X \rightarrow B$ such that
for each $s \in \Gamma^+$, $(\psi_s, \psi_0)$ is a Toeplitz
representation of $X_s$, and $\psi (xy) = \psi (x) \psi (y)$ for
$x, y \in X$.
\end{definition}

Every product system $X$ over $\Gamma^+$ has a \emph{Toeplitz
algebra} $\mathcal{T}_X$ which is generated by a universal
Toeplitz representation $i_X$ of $X$, and it is unique up to
isomorphism \cite[Proposition 2.8]{f}.

\section{Partial-isometric representations}\label{sectpireps}

Let $\Gamma$ be a totally ordered abelian group with positive cone
$\Gamma^+$. A \emph{partial-isometric representation} $V$ of
$\Gamma^+$ on a Hilbert space $H$ is a map $V$ of $\Gamma^+$ into
$B (H)$ such that $V_s$ is a partial isometry and $V_s V_t =
V_{s+t}$ for every $s,t$ in $\Gamma^+$. We denote by $C^*
(V(\Gamma^+))$ the $C^*$-algebra generated by the operators $V_s$.

\begin{example}
Since $V_n=(V_1)^n$, a partial-isometric representation $V$ of
$\N$ is determined by the single partial isometry $V_1$; a single
partial isometry $W$ generates a partial isometric representation
$V:n\mapsto W^n$ if and only $W$ is a power partial isometry. We
often implicitly acknowledge this by writing $V^n$ for $V_n$ when
the semigroup is $\N$.
\end{example}

The following property of partial-isometric representations will
be used repeatedly.

\begin{prop}\label{commprojs}
Suppose  $V$ is a partial-isometric representation of $\Gamma^+$
on $H$. Then each $V_s$ is a power partial isometry and $\{V_s^*
V_s, V_t V_t^*: s, t \in \Gamma^+\}$ is a commuting family of
projections.
\end{prop}

\begin{proof}
For each $s\in \Gamma^+$ and $n\in \mathbb{N}$, $V_s^n = V_{ns}$
is a partial isometry, so $V_s$ is a power partial isometry.
Because $V_sV_t=V_{s+t}$ is a partial isometry, each $V_s^*V_s$
commutes with each $V_tV_t^*$ by Proposition~\ref{tool}. To see
that the range projections commute, we let $s,t\in \Gamma^+$ and
compute:
\begin{align*}
V_s V_s^* V_t V_t^* &= \left\{ \begin{array}{ll}
V_s V_{s-t}^* V_t^* V_t V_t^* & \mbox{if $t \leq s$} \\
V_s V_s^* V_s V_{t-s} V_t^* & \mbox{if $s \leq t$}
\end{array}
\right. \\
&= \left\{ \begin{array}{ll}
V_s V_{s-t}^* V_t^* & \mbox{if $t \leq s$} \\
V_s V_{t-s} V_t^* & \mbox{if $s \leq t$}
\end{array}
\right. \\
&= \left\{ \begin{array}{ll}
V_s V_s^* & \mbox{if $t \leq s$} \\
V_t V_t^* & \mbox{if $s \leq t$}
\end{array}
\right. \\
&= V_{s \vee t} V_{s \vee t}^*,
\end{align*}
where $s \vee t$ denotes $\max \{s,t\}$. Since $s \vee t = t \vee
s$, this and the same calculation with $s$ and $t$ swapped show
that $V_s V_s^*$ commutes with $V_t V_t^*$. A similar argument
shows that $V_s^* V_s V_t^* V_t  = V_{s \vee t}^* V_{s \vee t}=
V_t^* V_t V_s^* V_s$.
\end{proof}

Every isometric representation $V:\Gamma^+\to B(H)$ is also a
partial-isometric representation, and so is the associated
coisometric representation $V^*:s\mapsto V_s^*$. Thus there are
natural representations $T$ and $T^*$ of $\Gamma^+$ by forward and
backward shifts on $\ell^2(\Gamma^+)$. In the remainder of this
section we discuss two partial-isometric representations by
truncated shifts, and the $C^*$-algebras they generate.

For $s \in \Gamma^+$, we consider two intervals
\[
[0,s):=\{t\in\Gamma^+:0\leq t<s\}\ \text{ and }\
[0,s]:=\{t\in\Gamma^+:0\leq t\leq s\}.
\]
For $t\in \Gamma^+$, there is a partial isometry $K_t^s$ on
$\ell^2 ([0,s))$ such that
\begin{equation}\label{defK} K^s_t (\varepsilon_r) = \left\{
\begin{array}{ll}
\varepsilon_{r+t} & \mbox{if $r+t\in [0,s)$} \\
0 & \mbox{otherwise,}
\end{array}
\right.
\end{equation}
and $K^s: \Gamma^+ \rightarrow B (\ell^2 ([0,s)))$ is a
partial-isometric representation of $\Gamma^+$ satisfying $K_t^s =
0$ for $t \geq s$. Similarly, there are partial isometries $J_t^s$
on $\ell^2([0,s])$ such that
\begin{equation}\label{defJ}
  J^s_t (\varepsilon_r) = \left\{
\begin{array}{ll}
\varepsilon_{r+t} & \mbox{if $r+t \in [0,s]$} \\
0 & \mbox{otherwise,}
\end{array}
\right.
\end{equation}
and then  $J^s: \Gamma^+ \rightarrow B (\ell^2 ([0,s]))$ is a
partial-isometric representation of $\Gamma^+$ satisfying $J_t^s =
0$ for $t
>  s$.

We analyse $C^* (K^s(\Gamma^+))$ first.

\begin{prop}\label{commidealCI}
There is an order ideal $I$ of $\Gamma$ such that
\begin{align}\label{idealI}
I^+ &= \{t \in \Gamma^+ : 0 \leq t \leq ns \mbox{ for some } n \in
\mathbb{N}\},
\end{align}
and then  $C^* (K^s (\Gamma^+))$ is Morita equivalent to the
commutator ideal $\mathcal{C}_I$ in $\mathcal{T} (I)$.
\end{prop}

Before we prove this Proposition we identify $C^* (K^s
(\Gamma^+))$ as a corner of the commutator ideal
$\mathcal{C}_{\Gamma}$.

\begin{lemma}\label{C*Ks}
For $s\in \Gamma^+$, $C^* (K^s (\Gamma^+))$ is isomorphic to $(1 -
T_sT_s^*) \mathcal{C}_{\Gamma} (1 - T_sT_s^*)$.
\end{lemma}

To prove Lemma~\ref{C*Ks}, we need the following standard fact.

\begin{lemma}\label{subspace}
Suppose $K$ is a closed subspace of a Hilbert space $H$ and $P$ is
the projection of $H$ onto $K$. Then $PTP\mapsto T|_K$ is an
isomorphism of $P B(H) P$ onto $B(K)$.
\end{lemma}

\begin{proof}[Proof of Lemma \ref{C*Ks}]
We view $\ell^2 ([0,s))$ as the closed subspace of $\ell^2
(\Gamma^+)$ spanned by $\{\varepsilon_t: t \in [0,s)\}$. Then $1 -
T_s T_s^*$ is the projection of $\ell^2 (\Gamma^+)$ onto $\ell^2
([0,s))$. We have
\[(1 - T_s T_s^*) T_t (1 - T_s T_s^*)
(\varepsilon_r) =
\begin{cases}
K^s_t (\varepsilon_r)&\text{if $r<s$}\\
0&\text{if $r\geq s$.}
\end{cases}\]
Thus the isomorphism of Lemma~\ref{subspace} identifies
$C^*(K^s(\Gamma^+))$ with the $C^*$-subalgebra
\[
D:=C^*\big(\{(1 - T_sT_s^*) T_t (1 - T_sT_s^*): t
\in\Gamma^+\}\big)
\]
of $\mathcal{T} (\Gamma)$. It therefore suffices to prove that $D
= (1 - T_sT_s^*) \mathcal{C}_{\Gamma} (1 -T_sT_s^*)$.

Since  $1 - T_sT_s^*=[T_s^*,T_s]$ belongs to
$\mathcal{C}_{\Gamma}$, and $\mathcal{C}_{\Gamma}$ is an ideal,
each $(1 - T_sT_s^*) T_t (1 - T_sT_s^*)$ belongs to $(1 -
T_sT_s^*) \mathcal{C}_{\Gamma} (1 - T_sT_s^*)$. Thus $D \subset (1
- T_sT_s^*) \mathcal{C}_{\Gamma} (1 - T_sT_s^*)$.

Before proving the reverse inclusion, we recall from
Lemma~\ref{comideal} that $\mathcal{C}_{\Gamma}$ is spanned by the
elements of the form $T_r (1 - T_u T_u^*) T_t^*$. Since
\begin{align*}
(1 &- T_s T_s^*) T_r (1 - T_u T_u^*) T_t^* (1 - T_s T_s^*) \\
&= \left\{ \begin{array}{ll} (1 - T_s T_s^*) T_{r-t} (T_t T_t^* -
T_{u+t} T_{u+t}^*) (1 - T_s T_s^*) & \mbox{if $0 \leq r-t
< s$} \\
(1 - T_s T_s^*) (T_r T_r^* - T_{u+r} T_{u+r}^*) T_{t-r}^* (1 - T_s
T_s^*) & \mbox{if $0 \leq
t-r < s$} \\
0 & \mbox{otherwise}
\end{array}
\right. \\
&= \left\{ \begin{array}{ll} (1 - T_s T_s^*) T_{r-t} (1 - T_s
T_s^*) (T_t T_t^* - T_{u+t} T_{u+t}^*) (1 - T_s T_s^*) &
\mbox{if $0 \leq r-t < s$} \\
(1 - T_s T_s^*) (T_r T_r^* - T_{u+r} T_{u+r}^*) (1 - T_s T_s^*)
T_{t-r}^* (1 - T_s T_s^*) &
\mbox{if $0 \leq t-r < s$} \\
0 & \mbox{otherwise,}
\end{array}
\right.
\end{align*}
it suffices to prove that $(1 - T_sT_s^*) T_tT_t^* (1 -
T_sT_s^*)\in D$  for every $t$. But since $T_t^*T_t=1$ and $1 -
T_sT_s^* \leq 1 - T_{t+s}T_{t+s}^*$, we calculate
\begin{align*}
(1 - T_sT_s^*) T_t (1 - T_sT_s^*) &= (1 - T_sT_s^*) (T_t -
T_tT_sT_s^*T_t^*T_t) \\
&= (1 - T_sT_s^*) (1 - T_{t+s}T_{t+s}^*) T_t \\
&= (1 - T_sT_s^*) T_t,
\end{align*}
and deduce that
\[(1 - T_sT_s^*) T_tT_t^* (1 - T_sT_s^*) =  (1 - T_sT_s^*) T_t ((1 - T_sT_s^*)
T_t)^*\] is in $D$. This proves the reverse inclusion, and hence
Lemma~\ref{C*Ks}.
\end{proof}

\begin{proof}[Proof of Proposition~\ref{commidealCI}]
Since $I^+$ is a subsemigroup of $\Gamma^+$ and $0 \leq t \leq r
\in I^+$ implies $t \in I^+$, $I:=I^+\cup (-I^+)$ is an order
ideal. For $t \in \Gamma^+ \setminus I^+$, we have $K_t^s = 0$.
Thus
\[C^* (K^s (\Gamma^+)) = C^*(\{K_t^s: t \in \Gamma^+\}) =
C^*(\{K_t^s: t \in I^+\}) = C^* (K^s (I^+)).\] By
Lemma~\ref{C*Ks}, $C^*(K^s (I^+))$ is isomorphic to the
$C^*$-subalgebra $(1 - V_s V_s^*) \mathcal{C}_I (1 - V_s V_s^*)$
of the commutator ideal $\mathcal{C}_I$, where, to avoid
eyestrain, we have written $V$ for $T^I$. But $(1 - V_s V_s^*)
\mathcal{C}_I (1 - V_s V_s^*)$ is Morita equivalent to the ideal
$\overline{\mathcal{C}_I (1 - V_s V_s^*) \mathcal{C}_I}$
\cite[Example~3.6]{rw}, so it suffices to prove that
$\mathcal{C}_I (1 - V_s V_s^*) \mathcal{C}_I $ is dense in
$\mathcal{C}_I$.

  Lemma~\ref{comideal} implies that
\begin{equation}\label{commutatorI}
\mathcal C_I = \overline{\newspan} \{V_r (1 - V_tV_t^*) V_u^*:
r,t,u \in I^+\}.
\end{equation}
Since $1 - V_s V_s^*$ is a projection in $\mathcal{C}_I$,  $1 -
V_s V_s^*$ belongs to $\mathcal{C}_I (1 - V_s V_s^*)
\mathcal{C}_I$. Now suppose $t \in I^+$, say $t\leq Ns$. An
induction argument shows that
\begin{align*}
1-V_{ns}V_{ns}^*&=(1-V_{(n-1)s}V_{(n-1)s}^*)+
(V_{(n-1)s}V_{(n-1)s}^*-V_{ns}V_{ns}^*)\\
&=(1-V_{(n-1)s}V_{(n-1)s}^*)+V_{(n-1)s}(1-V_sV_s^*)V_{(n-1)s}^*
\end{align*}
belongs to $\overline{\mathcal{C}_I (1 - V_s V_s^*)
\mathcal{C}_I}$ for all $n$, and hence so does $1-V_tV_t^*\leq
1-V_{Ns}V_{Ns}^*$. We deduce that $V_r (1 - V_t V_t^*) V_u^*$ is
in $\mathcal{C}_I (1 - V_s V_s^*) \mathcal{C}_I$ for every $r,t,u
\in I^+$, and (\ref{commutatorI}) implies that $\mathcal{C}_I (1 -
V_s V_s^*) \mathcal{C}_I$ is dense in $\mathcal{C}_I$.
\end{proof}

We now consider $C^* (J^s (\Gamma^+))$.

\begin{prop}\label{C*Js}
By identifying $\ell^2 ([0,s))$ with a closed subspace of $\ell^2
([0,s])$, we can view $C^* (K^s (\Gamma^+))$ as a $C^*$-subalgebra
of $C^* (J^s (\Gamma^+))$, and then
\begin{align}\label{JsKsK}
C^* (J^s (\Gamma^+)) &= C^* (K^s (\Gamma^+)) + \mathcal{K} (\ell^2
([0,s])).
\end{align}
\end{prop}

\begin{proof}
When we view $\ell^2 ([0,s))$ as
$\overline{\newspan}\{\varepsilon_r : r \in [0,s)\}\subset \ell^2
([0,s])$, $1 - J_s^s (J_s^s)^*$ is the projection of $\ell^2
([0,s])$ onto $\ell^2 ([0,s))$. By Lemma~\ref{subspace}, there is
an isomorphism of $B (\ell^2 ([0,s)))$ onto $(1 - J_s^s (J_s^s)^*)
B(\ell^2 ([0,s])) (1 - J_s^s (J_s^s)^*)$. Moreover, for $r < s$,
\begin{align*}
(1 - J_s^s (J_s^s)^*) J^s_t (1 - J_s^s (J_s^s)^*) (\varepsilon_r)
&= (1 - J_s^s (J_s^s)^*) J^s_t (\varepsilon_r) \\
&= \left\{ \begin{array}{ll}
(1 - J_s^s (J_s^s)^*) (\varepsilon_{r+t}) & \mbox{if $r+t\in [0,s]$} \\
0 & \mbox{otherwise}
\end{array}
\right. \\
&= \left\{ \begin{array}{ll}
\varepsilon_{r+t} & \mbox{if $r+t\in [0,s)$} \\
0 & \mbox{otherwise}
\end{array}
\right. \\
&= K^s_t (\varepsilon_r),
\end{align*}
and $(1 - J_s^s (J_s^s)^*) J^s_t (1 - J_s^s (J_s^s)^*)
(\varepsilon_s) = 0$; thus we can identify $C^* (K^s (\Gamma^+))$
with the $C^*$-subalgebra of $C^* (J^s (\Gamma^+))$ generated by
\[\{(1 - J_s^s (J_s^s)^*) J^s_t (1 - J_s^s (J_s^s)^*): t \in \Gamma^+\}.\]

We notice that the rank-one operator $\varepsilon_r \otimes
\overline{\varepsilon}_t = (J_{s-r}^s)^* J_s^s (J_s^s)^*
J_{s-t}^s$, and hence
\begin{equation}\label{idmatrixunits}
\mathcal{K} (\ell^2 ([0,s])) = \overline{\newspan} \{(J_{s-r}^s)^*
J_s^s (J_s^s)^* J_{s-t}^s: r,t \in [0,s]\} \subset C^* (J^s
(\Gamma^+)).
\end{equation}
Thus $C^* (K^s (\Gamma^+)) + \mathcal{K} (\ell^2 ([0,s])) \subset
C^* (J^s (\Gamma^+))$. Note that $C^* (K^s (\Gamma^+)) +
\mathcal{K} (\ell^2 ([0,s]))$ is a $C^*$-subalgebra of $C^* (J^s
(\Gamma^+))$ because $\mathcal{K} (\ell^2 ([0,s]))$ is an ideal.
On the other hand, we compute using the definitions of $J^s$ and
$K^s$ that
\[
(J^s_t - K^s_t) (\varepsilon_r) = \left\{
\begin{array}{ll}
(\varepsilon_s \otimes \overline{\varepsilon}_{s-t})
(\varepsilon_r)
& \mbox{if $t \leq s$} \\
0 & \mbox{otherwise}
\end{array}
\right.
\]
for every $r \in [0,s]$, and hence
\[ J^s_t = \left\{
\begin{array}{ll}
K^s_t + (\varepsilon_s \otimes \overline{\varepsilon}_{s-t}) &
\mbox{if $t \leq s$} \\
0 & \mbox{otherwise.}
\end{array}
\right. \] Thus $J^s_t$ belongs to $C^* (K^s (\Gamma^+)) +
\mathcal{K} (\ell^2 ([0,s]))$ for every $t \in \Gamma^+$, and $C^*
(J^s (\Gamma^+)) \subset C^* (K^s (\Gamma^+)) + \mathcal{K}
(\ell^2 ([0,s]))$.
\end{proof}

\begin{cor}\label{exactseqCJs}
Let $\mathcal{J} = C^* (K^s (\Gamma^+)) \cap \mathcal{K} (\ell^2
([0,s]))$. Then $\mathcal{J}$ is an ideal of $C^* (K^s
(\Gamma^+))$ and there is an exact sequence
\begin{equation}\label{Rs}
0\longrightarrow {\mathcal{K} (\ell^2 ([0,s]))}\longrightarrow
{C^* (J^s (\Gamma^+))} \mathop{\longrightarrow}\limits^{R_s} {C^*
(K^s (\Gamma^+))/ \mathcal{J}}\longrightarrow 0
\end{equation}
in which $R_s (J^s_t) = K^s_t + \mathcal{J}$.
\end{cor}

\begin{proof}
We trivially have that $\mathcal{J}$ is an ideal of $C^* (K^s
(\Gamma^+))$. By Lemma~\ref{C*Js}, $C^* (J^s (\Gamma^+)) = C^*
(K^s (\Gamma^+)) + \mathcal{K} (\ell^2 ([0,s]))$. Then $C^* (J^s
(\Gamma^+))/\mathcal{K} (\ell^2 ([0,s])) = C^* (K^s (\Gamma^+))/
\mathcal{J}$, and we have (\ref{Rs}).
\end{proof}

When $\Gamma$ is Archimedean, and hence isomorphic to a subgroup
of the additive group of real numbers, we can say more.

\begin{prop}\label{kerRs}
Suppose $\Gamma$ is Archimedean. Then either $\Gamma$ is singly
generated, in which case $C^* (J^s (\Gamma^+)) = \mathcal{K}
(\ell^2 ([0,s]))$, or $\Gamma$ is isomorphic to a dense subgroup
of $\mathbb{R}$, in which case we have an exact sequence
\begin{equation}\label{exseqCJs1}
0\longrightarrow {\mathcal{K} (\ell^2 ([0,s]))}\longrightarrow
{C^* (J^s (\Gamma^+))} \mathop{\longrightarrow}\limits^{R_s} {C^*
(K^s (\Gamma^+))}\longrightarrow 0
\end{equation}
such that $R_s (J^s_t) = K^s_t$.
\end{prop}

\begin{proof}
If $\Gamma$ is singly generated, then $\Gamma=\mathbb{Z}t$ for
some $t$, $s=nt$ for some $n\geq 0$, $J^s_t$ is the truncated
shift on $\mathbb{C}^{n+1}\cong\ell^2([0,s])$, and
$C^*(J^s(\Gamma^+))=C^*(J^s_t)$ is all of
$\mathcal{K}(\ell^2([0,s])=B(\ell^2([0,s])$ by
Example~\ref{ppisometryJk}. Now suppose $\Gamma$ is dense in
$\mathbb{R}$. Then $\Gamma$ is simple, and the order ideal $I$ of
Proposition~\ref{commidealCI} is all of $\Gamma$; since
$\mathcal{C}_{\Gamma}$ is simple, Proposition~\ref{commidealCI}
implies that $C^*(K^s)$ is simple too.  Thus
$\mathcal{J}:=C^*(K^s)\cap \mathcal{K}(\ell^2([0,s])$ is either
$0$ or $C^*(K^s)$. But $\Gamma$ is dense,
$\ell^2([t,s])=K^s_t(\ell^2([0,s]))$ is infinite-dimensional for
all $t<s$, and hence $K^s_t$ is not compact whenever $t<s$. Thus
$\mathcal{J}=0$, and Corollary~\ref{exactseqCJs} gives
  (\ref{exseqCJs1}).
\end{proof}

\section{Partial-isometric crossed products}\label{sect2}

\subsection{Covariant partial-isometric representations}\label{covpirep}

We consider a dynamical system $(A, \Gamma^+, \alpha)$ consisting
of an action $\alpha$ of $\Gamma^+$ by endomorphisms of a
$C^*$-algebra $A$ such that $\alpha_0 = \id_A$. We assume that
each $\alpha_s$ is extendible, and hence extends to a strictly
continuous endomorphism $\overline{\alpha}_s$ of the multiplier
algebra $M(A)$.

\begin{definition}
A \emph{covariant partial-isometric representation} of $(A,
\Gamma^+, \alpha)$ is a pair $(\pi, V)$ consisting of a
non-degenerate representation $\pi$ of $A$ on a Hilbert space $H$
and a partial-isometric representation $V$ of $\Gamma^+$ on $H$
such that
\begin{align}\label{covrep}
\pi (\alpha_s (a)) = V_s \pi(a) V_s^* \ \mbox{ and }\  V_s^* V_s
\pi (a) = \pi (a) V_s^* V_s \ \mbox{ for $s \in \Gamma^+$ and $a
\in A$.}
\end{align}
We can also talk about  covariant partial-isometric
representations of $(A, \Gamma^+, \alpha )$ in a $C^*$-algebra $B$
or a multiplier algebra $M(B)$.
\end{definition}

\begin{lemma}\label{covmultiplier}
If $(\pi, V)$ is a covariant partial-isometric representation of
$(A, \Gamma^+, \alpha)$, then $(\overline{\pi}, V)$ is a covariant
partial-isometric representation of $(M(A), \Gamma^+,
\overline{\alpha})$.
\end{lemma}

\begin{proof}
Let $\{a_{i}\}$ be an approximate identity for $A$.  Then
\[
\pi (\alpha_s (m a_{i})) = V_s \pi (m a_{i}) V_s^* = V_s
\overline{\pi} (m) \pi (a_{i}) V_s^*
\]
converges strongly to $V_s \overline{\pi} (m) V_s^*$, because
$\pi$ is non-degenerate. On the other hand, since $\alpha_s
(a_{i})$ converges strictly to $\overline{\alpha}_s (1)$,
\[
\pi (\alpha_s (m a_{i})) = \overline{\pi} (\overline{\alpha}_s
(m)) \pi (\alpha_s (a_{i}))
\]
converges strongly to $\overline{\pi} (\overline{\alpha}_s (m))
\pi (\alpha_s (1)) = \overline{\pi} (\overline{\alpha}_s (m))$.
Thus $\overline{\pi} (\overline{\alpha}_s (m)) = V_s
\overline{\pi} (m) V_s^*$. A similar argument shows that $V_s^*
V_s \overline{\pi} (m) = \overline{\pi} (m) V_s^* V_s$.
\end{proof}

Next, we give an alternative version of the covariance relation
(\ref{covrep}).

\begin{lemma}\label{equivalence}
Let $A$ be a C*-algebra and $\alpha$ be an extendible endomorphism
of $A$. Suppose $\pi$ is a non-degenerate homomorphism of $A$ into
a multiplier algebra $M(B)$ and $v$ is a partial isometry in
$M(B)$. Then
\begin{equation}\label{1stcov}
\pi (\alpha(a)) v = v \pi (a) \mbox{ for every $a \in A$}\mbox{
and } vv^* = \overline{\pi} (\overline{\alpha} (1))
\end{equation} if and only if
\begin{equation}\label{2ndcov}
\pi (\alpha(a)) = v \pi(a) v^* \mbox{ and } v^*v \pi(a) = \pi (a)
v^*v\ \mbox{for every $a \in A$.}
\end{equation}
\end{lemma}

\begin{proof}
Suppose (\ref{1stcov}) holds.  Then
\begin{gather*}
\pi (\alpha (a)) = \pi (\alpha (a)) \overline{\pi}
(\overline{\alpha} (1)) = \pi (\alpha (a)) vv^* = v
\pi (a)v^*,\text{ and}\\
v^*v \pi (a) = v^* \pi (\alpha (a))v = (\pi (\alpha (a^*))v)^* v =
(v \pi (a^*))^* v = \pi(a) v^*v.
\end{gather*}
Now  suppose (\ref{2ndcov}) holds, and $a\in A$. Then
\[\pi (\alpha (a)) v = v \pi (a)
v^*v = vv^*v \pi (a) = v \pi (a),\] and Lemma~\ref{covmultiplier}
gives $\overline{\pi} (\overline{\alpha} (1)) = v \overline{\pi}
(1) v^* = vv^*$.
\end{proof}

\begin{cor}\label{altcov}
Let $\pi$ be a non-degenerate representation of $A$ on a Hilbert
space $H$ and $V$ be a partial-isometric representation of
$\Gamma^+$ on $H$. Then $(\pi , V)$ is a covariant
partial-isometric representation of $(A, \Gamma^+, \alpha)$ if and
only if
\[\pi (\alpha_s (a)) V_s = V_s \pi (a)\ \mbox{ and }\ \overline{\pi}
(\overline{\alpha}_s (1)) = V_s V_s^*\ \mbox{ for $a\in A$, $s\in
\Gamma^+$.}\]
\end{cor}

\begin{cor}
If $(\pi, V)$ is a covariant partial-isometric representation of
$(A, \Gamma^+, \alpha)$ on $H$, then $V_0 = 1$.
\end{cor}

\begin{proof}
Using Proposition~\ref{commprojs}, we calculate
\begin{align*}
V_0 = V_{0+0} V_{0+0}^* V_0 = V_0 (V_0 V_0^*)( V_0^* V_0)= V_0
(V_0^* V_0)( V_0 V_0^*)= (V_0 V_0^*) V_{0+0} V_0^* = V_0 V_0^*.
\end{align*}
Since we assume that $\alpha_0$ is the identity endomorphism, and
$\pi$ is non-degenerate, we get $V_0 V_0^* = \overline{\pi}
(\overline{\alpha}_0 (1)) = 1$. Thus $V_0 = 1$.
\end{proof}

\begin{example}\label{repwithAfaithful}
Suppose $\alpha:\Gamma^+\to\End A$, and  $\pi_0$ is a
non-degenerate representation of $A$ on $H$. Define $\pi:A\to
B(\ell^2 (\Gamma^+, H))$ by
\[(\pi (a) \zeta) (r) = \pi_0 (\alpha_r (a)) (\zeta (r)).\]
Then $\pi$ is a representation which is non-degenerate on
\[\mathcal{H} = \overline{\newspan} \{\zeta\in \ell^2 (\Gamma^+,
H) : \zeta (r) \in \overline{\pi}_0 (\overline{\alpha}_r (1))H
\mbox{ for all } r \in \Gamma^+\}.\] For $s \in \Gamma^+$, define
$V_s$ on $\mathcal{H}$ by $V_s \zeta (r) = \zeta (r+s)$: for
$\zeta \in \mathcal{H}$, $V_s \zeta (r)$ is in $\overline{\pi}_0
(\overline{\alpha}_{r+s} (1))H \subset \overline{\pi}_0
(\overline{\alpha}_r (1))H$, and hence $V_s \zeta \in
\mathcal{H}$. Let $\zeta \in (\ker V_s)^{\perp}$. Then $\zeta (r)
=0$ for every $r < s$, and hence
\[\|V_s \zeta\|^2 =
\sum_{r \in \Gamma^+} \|\zeta (r+s)\|^2 =\sum_{r \in \Gamma^+}
\|\zeta (r)\|^2  = \|\zeta\|^2.\] Thus $V_s$ is a partial
isometry. Indeed, $V:s\mapsto V_s$ is a partial-isometric
representation of $\Gamma^+$ on $\mathcal{H}$, because \[V_{s+t}
\zeta (r) = \zeta (r+s+t) = V_t \zeta (r+s) = V_s V_t \zeta (r)
\text{ for $\zeta \in \mathcal{H}$ and $s,t,r \in \Gamma^+$.}
\]

We claim that $(\pi|_{\mathcal{H}}, V)$ is a covariant
partial-isometric representation of $(A, \Gamma^+, \alpha)$. For
$a \in A$, $\zeta \in \mathcal{H}$ and $r,s \in \Gamma^+$, we have
\begin{align*}
\big(\pi (\alpha_s (a)) (V_s \zeta)\big) (r) &= \pi_0 (\alpha_r
(\alpha_s (a)))(V_s \zeta (r))
= \pi_0 (\alpha_{r+s} (a)) (\zeta (r+s)) \\
&= (\pi (a) \zeta) (r+s) = (V_s \pi (a) \zeta) (r),
\end{align*}
and hence $\pi (\alpha_s (a)) V_s = V_s \pi (a)$. Since $V_s
V_s^*$ is the projection on
\[\range V_s = \{\zeta \in \mathcal{H}: \zeta (r) \in \overline{\pi}_0
(\overline{\alpha}_{r+s} (1)) H \mbox{ for all } r \in
\Gamma^+\},\] we get
\begin{align*}
(\overline{\pi} (\overline{\alpha}_s (1)) \zeta) (r) =
\overline{\pi}_0 (\overline{\alpha}_r (\overline{\alpha}_s (1)))
(\zeta (r)) = \overline{\pi}_0 (\overline{\alpha}_{r+s} (1))
(\zeta (r))= (V_s V_s^* \zeta) (r)
\end{align*}
for every $\zeta \in \mathcal{H}$ and $r \in \Gamma^+$.
Corollary~\ref{altcov} implies that $(\pi, V)$ is covariant, as
claimed.

Notice that if $\pi_0$ is fathful, so is $\pi$. Thus every system
$(A,\Gamma^+,\alpha)$ admits covariant partial-isometric
represntations $(\pi,V)$ with $\pi$ faithful.
\end{example}

\subsection{Crossed products and the Toeplitz algebras of Hilbert bimodules}

Suppose $(A, \Gamma^+, \alpha)$ is a dynamical system as in
\S\ref{covpirep}. Following Fowler \cite[\S3]{f}, we give each
$X_s:= \{s\} \times \overline{\alpha}_s (1) A$ the structure of a
Hilbert bimodule over $A$ via
\begin{align*}
(s, x) \cdot a := (s, xa),\quad \langle (s, x), (s, y) \rangle_A
:= x^*y\quad\text{and}\quad a\cdot (s,x) = (s, \alpha_s (a) x),
\end{align*}
and define a multiplication on $X = \bigsqcup\{X_s : s \in
\Gamma^+\}$ by
\[(s, x) (t, y) = (s+t, \alpha_t (x) y)\ \mbox{ for $x \in
\overline{\alpha}_s (1) A$ and $y \in \overline{\alpha}_t (1) A$.}
\]
Then $X$ is a product system of Hilbert bimodules over $\Gamma^+$
\cite[Lemma 3.2]{f}. The Toeplitz algebra $(\mathcal{T}_X,i_X)$ of
\cite[Proposition 2.8]{f} is universal for covariant
partial-isometric representations:

\begin{prop}\label{crossedproduct}
Let $\alpha:\Gamma^+\to \End A$ be an action by extendible
endomorphisms. Then there is a covariant partial-isometric
representation $(i_A, i_{\Gamma^+})$ of $(A, \Gamma^+,\alpha)$ in
$\mathcal{T}_X$ such that $i_A$ is injective and
\smallskip

\textnormal{(a)} for every covariant partial-isometric
representation $(\pi, V)$ of $(A, \Gamma^+, \alpha)$ on $H$, there
is a non-degenerate representation $\pi \times V$ of
$\mathcal{T}_X$ on $H$ such that $(\pi \times V) \circ i_A = \pi$
and $(\overline{\pi \times V}) \circ i_{\Gamma^+} = V$; and
\smallskip

\textnormal{(b)} $\mathcal{T}_X$ is generated by $i_A (A) \cup
i_{\Gamma^+} (\Gamma^+)$; indeed, we have
\begin{equation}\label{spanningset}
\mathcal{T}_X=\overline{\newspan}\{i_{\Gamma^+}(s)^*
i_A(a)i_{\Gamma^+}(t): a\in A,\ s,t\in\Gamma^+\}.
\end{equation}

\noindent If $(j_A,j_{\Gamma^+})$ is a covariant partial-isometric
representation of $(A,\Gamma^+,\alpha)$ in a $C^*$-algebra $B$
with properties \textnormal{(a)} and \textnormal{(b)}, then there
is an isomorphism of $B$ onto $\mathcal{T}_X$ which carries
$(j_A,j_{\Gamma^+})$ into $(i_A, i_{\Gamma^+})$.
\end{prop}

\begin{proof}
By \cite[Proposition~3.4]{f}, there is a partial-isometric
representation $(i_A, i_{\Gamma^+})$ in $\mathcal{T}_X$ such that
$i_A (a) = i_X (0, a)$ and $i_{\Gamma^+} (s) = \lim i_X (s,
\alpha_s (a_i))^*$, where $\{a_i\}$ is an approximate identity for
$A$. Since \cite[Proposition~3.4]{f} also says that we can recover
$i_X$ via $i_X(s,a)=i_{\Gamma^+}(s)^*i_A(a)$, the elements $i_A
(a)$ and $i_{\Gamma^+} (s)$ generate. To verify
(\ref{spanningset}), just check that the right-hand side is closed
under multiplication, and hence is a $C^*$-algebra containing the
generators.

Let $(\pi , V)$ be a partial-isometric representation of $(A,
\Gamma^+, \alpha)$ on $H$. Then
  \cite[Proposition~3.4]{f} gives a Toeplitz representation $\psi$ of
$X$ such that $\psi (s,x) = V_s^* \overline{\pi} (x)$, and by
\cite[Proposition 2.8(a)]{f} there is a representation $\psi_*$ of
$\mathcal{T}_X$ with $\psi_*\circ i_X=\psi$. We take $\pi \times V
= \psi_*$. Then $\pi \times V (i_A (a)) = \pi \times V (i_X (0,a))
= \psi (0,a) = \pi (a)$, and
\begin{align*}
\pi \times V (i_{\Gamma^+} (s)) &= \pi \times V (\lim i_X (s,
\alpha_s (a_i))^*) \\
&= \lim \pi \times V (i_X (s, \alpha_s (a_i))^*) \\
&= \lim \psi_*(i_X (s, \alpha_s (a_i))^*) \\
&= \lim \psi (s, \alpha_s (a_i))^* \\
&= V_s.
\end{align*}
Thus $(i_A, i_{\Gamma^+})$ satisfies (a).
Example~\ref{repwithAfaithful} shows that there are covariant
representations $(\pi,V)$ with $\pi$ faithful, and then the
equation $\pi=(\pi\times V)\circ i_A$ shows that $i_A$ is
injective. The uniqueness follows by a standard argument using the
two universal properties.
\end{proof}

We call $\mathcal{T}_X$ the \emph{partial-isometric crossed
product} of $(A, \Gamma^+, \alpha)$, and denote it $A
\times_{\alpha} \Gamma^+$. Because it is the Toeplitz algebra of a
product system, we can apply Fowler's results, and in particular
\cite[Corollary 9.4]{f}. This gives:

\begin{thm}\label{faithful}
Suppose $(\pi, V)$ is a covariant partial-isometric representation
of the system $(A, \Gamma^+, \alpha)$ on $H$. Then  $\pi \times V$
is faithful on $A\times_{\alpha}\Gamma^+$ if and only if $\pi$
acts faithfully on $(V_s^*H)^{\perp}$ for every $s$ in $\Gamma^+
\setminus \{0\}$.
\end{thm}

\section{The crossed product of $(B_{\Gamma^+}, \Gamma^+,
\tau)$}\label{sect3}

In this section, we analyse the crossed product
$B_{\Gamma^+}\times_\tau\Gamma^+$ by the action $\tau$ of
$\Gamma^+$ by right translation on $B_{\Gamma^+}\subset
\ell^\infty(\Gamma^+)$, which by \cite[Proposition~9.6]{f} is
universal for partial-isometric representations of $\Gamma^+$.
Theorem~\ref{Bgammatau} concerns subgroups of $\mathbb{R}$, but
until then $\Gamma$ can be any totally ordered abelian group.

We begin by analysing some related crossed products associated to
intervals in $\Gamma^+$.  Let $s \in\Gamma^+$, and let $I$ stand
for one of $[0,s)$ or $[0,s]$. For $t \in I$, let $1_t^I$ be the
characteristic function of $I\cap [t,\infty)$. Since $(1^I_t)^* =
1^I_t$ and $1^I_t 1^I_u = 1^I_{\max \{t,u\}}$, $B_I =
\overline{\newspan} \{1^I_t: t \in I\}$ is a $C^*$-subalgebra of
$\ell^{\infty} (I)$. We denote by $\tau^I$ the action of the
semigroup $\Gamma^+$ by translation on $\ell^{\infty} (I)$:
\[ \tau_r^I (f) (t) = \left\{
\begin{array}{ll}
f (t-r) & \mbox{if $t-r \in I$} \\
0 & \mbox{otherwise.}
\end{array}
\right. \] Then $\tau_r^I (1_t^I)$ is $1_{r+t}^I$ if $r+t \in I$,
and $0$ otherwise, so $B_I$ is invariant under $\tau^I$. We thus
obtain a dynamical system $(B_I, \Gamma^+, \tau^I)$. Since
$\tau^I_r=0$ when $r$ is not in $I$,  $B_I \times_{\tau^I}
\Gamma^+$ for $I=[0,s]$ or $[0,s)$ is quite different in nature
from $B_{\Gamma^+} \times_{\tau} \Gamma^+$. Nevertheless, the
crossed products $B_I \times_{\tau^I} \Gamma^+$ play an important
role in our structure theorem for $B_{\Gamma^+} \times_{\tau}
\Gamma^+$, and hence we analyse $B_I \times_{\tau^I} \Gamma^+$
first. We begin by describing their universal property.

\begin{prop}\label{univbtaugamma}
The crossed product $B_I \times_{\tau^I} \Gamma^+$ is the
universal $C^*$-algebra generated by a partial-isometric
representation $V$ of $\Gamma^+$ such that $V_r = 0$ for $r \notin
I$.
\end{prop}

We let $(i_{B_I}, i_{\Gamma^+}^I)$ denote the universal covariant
partial-isometric representation which generates $B_I
\times_{\tau^I} \Gamma^+$. The covariance relation $i_{B_I}
(1_r^I) = i_{\Gamma^+}^I (r) i_{\Gamma^+}^I (r)^*$ implies that
the partial-isometric representation $i_{\Gamma^+}^I$ itself
generates $B_I \times_{\tau^I} \Gamma^+$, and also that
$i_{\Gamma^+}^I (r) = 0$ for $r \notin I$. So for any covariant
partial-isometric representation $(\pi, V)$ of $(B_I, \Gamma^+,
\tau^I)$ we have $V_r=0$ when $r$ is not in $I$. It remains to
prove that if $V$ is a partial-isometric representation of
$\Gamma^+$ such that $V_r = 0$ for $r \notin I$, then there is a
representation $\pi_V^I$ of $B_I$ such that $(\pi_V^I,V)$ is
covariant, because then $V=(\pi_V^I\times V)\circ i_{\Gamma^+}^I$
factors through $i_{\Gamma^+}^I$. Thus the Proposition will follow
from  Lemma~\ref{excovrep} below.

For the proof of Lemma~\ref{excovrep}, we need the following
variant of \cite[Proposition~2.2]{alnr} and
\cite[Proposition~1.3]{lr}, which can be proved by making minor
modifications to the proof of \cite[Proposition~2.2]{alnr}.

\begin{lemma}\label{projtorep}
Suppose $\{P_r:r\in I\}$ is a family of projections on $H$ such
that $P_r\geq P_t$ when $r\leq t$. Then there is a representation
$\pi_P$ of $B_I$ on $H$ such that $\pi_P(1^I_r)=P_r$ for $r\in I$,
and $\pi_P$ is faithful if and only if $P_r\not=P_t$ when
$r\not=t$.
\end{lemma}

\begin{lemma}\label{excovrep}
Let $V$ be a partial-isometric representation of $\Gamma^+$ on $H$
such that $V_r = 0$ when $r \notin I$. Then there is a
representation $\pi_V^I$ of $B_I$ on $H$ such that $(\pi_V^I, V)$
is a covariant partial-isometric representation of the dynamical
system $(B_I, \Gamma^+, \tau^I)$.
\end{lemma}

\begin{proof}
First we prove that there is a representation $\pi_V^I$ of $B_I$
such that $\pi_V^I (1_r^I) = V_r V_r^*$. For $t > r$, we have
\begin{align*}
V_r V_r^* - V_t V_t^* = V_r (1 - V_{t-r} V_{t-r}^*) V_s^* = ((1 -
V_{t-r} V_{t-r}^*) V_r^*)^* (1 - V_{t-r} V_{t-r}^*) V_r^*,
\end{align*}
so $V_r V_r^* \geq V_t V_t^*$. Lemma~\ref{projtorep} now gives the
required representation.

To see that $(\pi_V^I, V)$ is covariant, it suffices by
Corollary~\ref{altcov} to show that $\pi_V^I (\tau_r^I (a)) V_r =
V_r \pi_V^I (a)$ for $r \in \Gamma^+$ and $a \in B_I$. By
continuity, we need only do this for $a=1^I_t$. For $r \notin I$
both sides of the equation are zero. For $r\in I$, we calculate
using Proposition~\ref{commprojs}: if $t+r\in I$, we have
\begin{align*}
\pi_V^I (\tau_r^I (1_t^I)) V_r = \pi_V^I (1_{t+r}^I) V_r= V_r (V_t
V_t^*) (V_r^* V_r) =(V_r V_r^* V_r) (V_t V_t^*) = V_r \pi_V^I
(1_t^I),
\end{align*}
and if $t+r \notin I$, we have \[\pi_V^I (\tau_r^I (1_t^I)) V_r =
\pi_V^I (0) V_r = 0 = V_{t+r} V_t^* = V_r V_t V_t^* = V_r \pi_V^I
(1_t^I).\] Thus $(\pi_V^I, V)$ is covariant.
\end{proof}

This completes the proof of Proposition~\ref{univbtaugamma}.

\begin{prop}\label{faithgammatau}
Let $V$ be a partial-isometric representation of $\Gamma^+$ such
that $V_r=0$ for $r\notin I$. Then the representation $\pi_V^I
\times V$ of $B_I \times_{\tau^I} \Gamma^+$ is faithful if and
only if \[(1 - V_r^* V_r) (V_u V_u^* - V_t V_t^*) \neq 0
\quad\text{for every $r > 0$, $u,t \in I$ and $u < t$.}
\]
\end{prop}

\begin{proof}
By Theorem~\ref{faithful}, $\pi_V^I \times V$ is faithful if and
only if $\pi_V^I|_{\range(1 - V_r^* V_r)}$ is faithful for every
$r > 0$. Let $r > 0$ and set $P_u:= (1 - V_r^* V_r) V_u V_u^*$ for
$u \in I$, which is a projection by Proposition~\ref{commprojs}.
The same Proposition implies that for $u \leq t$,
\begin{align*}
P_u - P_t &= (1 - V_r^* V_r) V_u V_u^* (1 - V_r^* V_r) - (1 -
V_r^* V_r) V_t V_t^* (1 - V_r^*
V_r) \\
&= (1 - V_r^* V_r) (V_u V_u^* - V_u V_{t-u} V_{t-u}^* V_u^*) (1 - V_r^* V_r) \\
&= (1 - V_r^* V_r) V_u (1 - V_{t-u} V_{t-u}^*) V_u^* (1 - V_r^* V_r) \\
&\geq 0.
\end{align*}
Thus by Lemma~\ref{projtorep}, there is a representation $\pi_P$
of $B_I$ on $(1-V_rV_r^*)H$ such that $\pi_P (1_u^I) = P_u$. Let
$h \in (1 - V_r^* V_r)(H)$. Then
\[
\pi_P (1_u^I) h = (1 - V_r^* V_r) V_u V_u^* h = V_u V_u^* (1 -
V_r^* V_r)h = V_u V_u^* h = \pi_V^I (1_u^I) h,
\]
and since the $1_u^I$ generate $B_I$, we deduce that $\pi_P =
\pi_V^I |_{(1 - V_r^* V_r)(H)}$. The Proposition therefore follows
from the second part of Lemma~\ref{projtorep}.
\end{proof}

\begin{cor}\label{JsIfaithful}
Let $J^r$ be the partial-isometric representation of $\Gamma^+$
satisfying \textnormal{(\ref{defJ})}. Then the representation
$\bigoplus_{r \in I} (\pi_{J^r}^I \times J^r)$ of
$B_I\times_{\tau^I}\Gamma^+$ on $\bigoplus_{r\in I} \ell^2
([0,r])$ is faithful.
\end{cor}

\begin{proof}
By Proposition~\ref{faithgammatau}, it is faithful if and only if
\begin{align}\label{forfaith}
\textstyle{\bigoplus_{r \in I}} \big(1 - (J^r_v)^* J^r_v\big)
\big(J^r_u (J^r_u)^* - J^r_t (J^r_t)^*\big) \neq 0
\end{align}
for every $v>0$ and $0\leq u < t\in I$. But the summand $\big(1 -
(J^u_v)^* J^u_v\big) \big(J^u_u (J^u_u)^* - J^u_t (J^u_t)^*\big)$
is nonzero, so  (\ref{forfaith}) holds.
\end{proof}

The crossed products $B_I\times_{\tau^I}\Gamma^+$ are important
because they arise as quotients of $B_{\Gamma^+}\times_\tau
\Gamma^+$. Each generating semigroup $i_{\Gamma^+}^I:\Gamma^+\to
B_I\times_{\tau^I}\Gamma^+$ is a partial-isometric representation,
and the universal property of $B_{\Gamma^+}\times_\tau \Gamma^+$
gives surjective homomorphisms $q_I:B_{\Gamma^+}\times_\tau
\Gamma^+\to B_I\times_{\tau^I}\Gamma^+$ such that
$q_I(i_{\Gamma^+}(r))=i_{\Gamma^+}^I(r)$.

For $s>0$, we write $q_s:=q_{[0,s]}$ and $q_s^-:=q_{[0,s)}$; for
$s=0$, we have only $q_0:= q_{[0,0]}$. Note that
\[\ker q_r^- \subset \ker q_s \subset \ker q_s^- \subset \ker q_t
\subset \ker q_0 \mbox{  for } t<s<r.\] Our structure theorem for
$B_{\Gamma^+} \times_{\tau} \Gamma^+$ involves these ideals and
the natural homomorphisms $\phi_T:=\pi_T^{\Gamma^+} \times
T:B_{\Gamma^+} \times_{\tau} \Gamma^+\to \mathcal{T}(\Gamma)$
associated to the Toeplitz representation $T:\Gamma^+\to
\mathcal{T}(\Gamma)$ and $\phi_{T^*}:=\pi_{T^*}^{\Gamma^+} \times
T^*:B_{\Gamma^+} \times_{\tau} \Gamma^+\to \mathcal{T}(\Gamma)$
associated to its adjoint $T^*:r\mapsto T^*_r$; both $\phi_T$ and
$\phi_{T^*}$ are surjective because the $T_r$ generate
$\mathcal{T}(\Gamma)$. The theorem also involves the homomorphisms
$\psi_T$ and $\psi_{T^*}$ of $\mathcal{T} (\Gamma)$ onto $C
(\hat{\Gamma})$ which carry $T_r$ to the evaluation maps
$\epsilon_r:\gamma\mapsto \gamma(r)$ and
$\epsilon_{-r}:\gamma\mapsto \gamma(-r)$.

\begin{thm}\label{Bgammatau}
Suppose $\Gamma$ is a subgroup of $\mathbb{R}$. Let
\[
\mathcal{I}=\big(\ker (\pi_T^{\Gamma^+} \times
T)\big)\cap\big(\ker (\pi_{T^*}^{\Gamma^+} \times T^*)\big).
\]
Then we have a commutative diagram
\begin{equation*}
\begin{diagram}\dgARROWLENGTH=0.8\dgARROWLENGTH
\node{} \node{0} \arrow{s} \node{0} \arrow{s} \node{0} \arrow{s}
\\
\node{0} \arrow{e} \node{\mathcal{I}} \arrow{e} \arrow{s}
\node{\ker (\pi_{T^*}^{\Gamma^+} \times T^*)} \arrow{e} \arrow{s}
\node{\mathcal{C}_{\Gamma}} \arrow{e} \arrow{s} \node{0}
\\
\node{0} \arrow{e} \node{\ker (\pi_T^{\Gamma^+} \times T)}
\arrow{e} \arrow{s} \node{B_{\Gamma^+} \times_{\tau} \Gamma^+}
\arrow{e,t}{\phi_T} \arrow{s,l}{\phi_{T^*}} \node{\mathcal{T}
(\Gamma)} \arrow{e} \arrow{s,l}{\psi_{T^*}} \node{0}
\\
\node{0} \arrow{e} \node{\mathcal{C}_{\Gamma}} \arrow{e} \arrow{s}
\node{\mathcal{T} (\Gamma)} \arrow{e,t}{\psi_{T}} \arrow{s}
\node{C (\hat{\Gamma})} \arrow{e} \arrow{s} \node{0}
\\
\node{} \node{0} \node{0} \node{0}
\end{diagram}
\end{equation*}
in which all the rows and columns are exact. For $s\in \Gamma^+$,
let $\mathcal{I}_s := \mathcal{I}\cap (\ker q_s)$,  and for $s>0$,
let $\mathcal{I}_s^- := \mathcal{I}\cap (\ker q_s^-)$. Then
\smallskip

\textnormal{(a)}  $\mathcal{I} / \mathcal{I}_0 \cong \mathbb{C}$;

\smallskip
\textnormal{(b)} $\mathcal{I}_s^- / \mathcal{I}_s \cong
\mathcal{K} (\ell^2 ([0,s]))$ for every $s>0$;
\smallskip

\textnormal{(c)} $\mathcal{I}_s^- = \bigcap_{t<s} \mathcal{I}_t$
for every $s>0$; and
\smallskip

\textnormal{(d)} $\mathcal{I}_s = \overline{\bigcup_{r > s}
\mathcal{I}_r^-}$ for every $s\in \Gamma^+$.
\end{thm}

The proof of this Theorem will occupy the rest of the section.

The right-hand exact sequence is due to Douglas \cite[Proposition
3]{d}. If $\Psi$ denotes the automorphism of $C (\hat{\Gamma})$
induced by the homeomorphism $\gamma\mapsto \gamma^{-1}$, then
$\psi_{T^*} = \Psi \circ \psi_T$, and hence the bottom sequence is
also exact. Since the middle sequences are exact by definition of
$\phi_T$ and $\phi_{T^*}$, we have the following diagram of exact
sequences:
\begin{equation*}
\begin{diagram}\dgARROWLENGTH=0.8\dgARROWLENGTH
\node{} \node{0} \arrow{s} \node{0} \arrow{s} \node{0} \arrow{s}
\\
\node{0} \arrow{e} \node{\mathcal{I}} \arrow{e} \arrow{s}
\node{\ker {\phi_{T^*}}}  \arrow{s} \node{\mathcal{C}_{\Gamma}}
\arrow{s}
\\
\node{0} \arrow{e} \node{\ker \phi_T} \arrow{e} \node{B_{\Gamma^+}
\times_{\tau} \Gamma^+} \arrow{e,t}{\phi_T}
\arrow{s,l}{\phi_{T^*}} \node{\mathcal{T} (\Gamma)} \arrow{e}
\arrow{s,l}{\psi_{T^*}} \node{0}
\\
\node{0} \arrow{e} \node{\mathcal{C}_{\Gamma}} \arrow{e}
\node{\mathcal{T} (\Gamma)} \arrow{e,t}{\psi_{T}} \arrow{s}
\node{C (\hat{\Gamma})} \arrow{e} \arrow{s} \node{0}
\\
\node{} \node{} \node{0} \node{0}
\end{diagram}
\end{equation*}
The top left-hand square commutes because all the maps are
inclusions, and the bottom right-hand square commutes because
\[
\psi_{T^*} \circ \phi_T (i_{\Gamma^+}
(t))=\psi_{T^*}(T_t)=\epsilon_{-t}=\epsilon_t^* =
\psi_{T}(T_t^*)=\psi_T\circ \phi_{T^*}(i_{\Gamma^+} (t))
\]
for every $t\in \Gamma^+$. The equation $\psi_{T^*} \circ
\phi_T=\psi_T\circ \phi_{T^*}$ also implies that $\phi_T$ maps
$\ker\phi_{T^*}$ into $\ker\psi_{T^*}=\mathcal{C}_{\Gamma}$, and
it maps $\ker\phi_{T^*}$ onto $\mathcal{C}_{\Gamma}$ because each
of the spanning elements $T_r(1 - T_vT_v^*) T_t^*$ in
Lemma~\ref{comideal} has the form $\phi_T(b)$ for $b=i_{\Gamma^+}
(r)(1 - i_{\Gamma^+} (v) i_{\Gamma^+} (v)^*) i_{\Gamma^+} (t)^*$
in $\ker\phi_{T^*}$. Since $\ker(\phi_T|_{\ker\phi_{T^*}})$ is by
definition $\mathcal{I}$, this gives exactness of the top row, and
exactness of the left-hand column follows similarly.

It remains to prove the assertions about the structure of
$\mathcal{I}$. Of these, (a) is easy: the homomorphism $q_0$ is
nonzero on the elements $(1 - i_{\Gamma^+} (u)^* i_{\Gamma^+} (u))
(1 - i_{\Gamma^+} (v) i_{\Gamma^+} (v)^*)$ of $\mathcal{I}$, and
has one-dimensional range. For the next two parts, we need a
lemma.

\begin{lemma}\label{kerqI}
For each interval $I$, $\ker q_I = \bigcap_{r \in I} \ker
(\pi_{J^r}^{\Gamma^+} \times J^r)$.
\end{lemma}

\begin{proof}
For $t \in \Gamma^+$, we have
\[\textstyle{\bigoplus_{r \in I} (\pi_{J^r}^{\Gamma^+} \times J^r)
(i_{\Gamma^+} (t)) = \bigoplus_{r \in I} J^r_t = \bigoplus_{r \in
I} (\pi_{J^r}^I \times J^r) \circ q_I (i_{\Gamma^+} (t)),}
\]
Since $\bigoplus_{r \in I} (\pi_{J^r}^I \times J^r)$ is faithful
on $B_I\times_{\tau^I}\Gamma^+$ by Corollary~\ref{JsIfaithful}, it
follows that $\ker q_I = \ker \bigoplus_{r \in I}
(\pi_{J^r}^{\Gamma^+} \times J^r)$.
\end{proof}

We now prove the remaining parts of Theorem~\ref{Bgammatau}. It is
convenient to do (c) first.

\begin{proof}[Proof of \textnormal{(c)}.]
 From Lemma~\ref{kerqI}, we have
\[
\textstyle{\bigcap_{t < s} \ker q_t = \bigcap_{t < s}
\big(\bigcap_{r \leq t} \ker (\pi_{J^r}^{\Gamma^+} \times
J^r)\big) = \bigcap_{r < s} \ker (\pi_{J^r}^{\Gamma^+} \times J^r)
= \ker q_s^-},
\]
and intersecting with $\mathcal{I}$ gives (c).
\end{proof}

\begin{proof}[Proof of \textnormal{(b)}.]
We shall prove that $\pi_{J^s}^{\Gamma^+} \times J^s$ is a
surjection of $\ker q_s^-$ onto $\mathcal{K} (\ell^2 ([0,s]))$
with kernel $\ker q_s$.  From (\ref{idmatrixunits}) we see that
$\mathcal{K} (\ell^2 ([0,s]))$ is spanned by the elements
\[
(J^s_{s-t})^* J^s_s (J^s_s)^* J^s_{s-r} = \pi_{J^s}^{\Gamma^+}
\times J^s (i_{\Gamma^+} (s-t)^* i_{\Gamma^+} (s) i_{\Gamma^+}
(s)^* i_{\Gamma^+} (s-r))
\]
of $\pi_{J^s}^{\Gamma^+} \times J^s(\ker q_s^-)$, so
$\pi_{J^s}^{\Gamma^+} \times J^s(\ker
q_s^-)\supset\mathcal{K}(\ell^2 ([0,s]))$. We next show the
reverse inequality.

If $\Gamma$ is singly generated, then $C^* (J^s (\Gamma^+)) =
\mathcal{K} (\ell^2 ([0,s]))$ by Example~\ref{ppisometryJk}, so
suppose $\Gamma$ is not singly generated. Then by
Proposition~\ref{kerRs} there is a homomorphism $R_s$ of $C^* (J^s
(\Gamma^+))$ onto $C^* (K^s (\Gamma^+))$ such that $R_s (J^s_t) =
K^s_t$ and $\ker R_s = \mathcal{K} (\ell^2 ([0,s]))$. But then we
can verify on generators that
\[
(\pi_{K^s}^{[0,s)} \times K^s) \circ q_s^-= R_s \circ
(\pi_{J^s}^{\Gamma^+} \times J^s),
\]
and this implies that $\pi_{J^s}^{\Gamma^+} \times J^s (\ker
q_s^-) \subset \ker R_s = \mathcal{K} (\ell^2 ([0,s]))$, as
claimed.

 From two applications of Lemma~\ref{kerqI} and (c), we obtain
\begin{align*}
\ker (\pi_{J^s}^{\Gamma^+} \times J^s) \cap (\ker q_s^-) &= \ker
(\pi_{J^s}^{\Gamma^+} \times J^s) \cap
\big(\textstyle{\bigcap_{t<s}} \ker(\pi_{J^t}^{\Gamma^+} \times
J^t)\big) \\
&= \textstyle{\bigcap_{t \leq s}} \ker(\pi_{J^t}^{\Gamma^+} \times J^t)\\
&= \ker q_s,
\end{align*}
and hence $\ker (\pi_{J^s}^{\Gamma^+} \times J^s|_{\ker q_s^-}) =
\ker q_s$.

We have now proved that $\pi_{J^s}^{\Gamma^+} \times J^s$ induces
an isomorphism of $(\ker q_s^-)/(\ker q_s)$ onto $\mathcal{K}
(\ell^2 ([0,s]))$. However, for any $r>0$ and $t>s$ the element
\[
\big(i_{\Gamma^+}(s)i_{\Gamma^+}(s)^*-i_{\Gamma^+}(t)i_{\Gamma^+}(t)^*\big)
\big(1-i_{\Gamma^+}(r)^*i_{\Gamma^+}(r)\big)
\]
belongs to $\mathcal{I}_s^-=(\ker q_s^-)\cap \mathcal{I}$ but not
to $\ker q_s$, and hence has nonzero image in $\mathcal{K} (\ell^2
([0,s]))$. Since the image of $\mathcal{I}_s^-$ is an ideal in
$\mathcal{K} (\ell^2 ([0,s]))$ and $\mathcal{K} (\ell^2 ([0,s]))$
is simple, this image must be all of $\mathcal{K} (\ell^2
([0,s]))$. So $\pi_{J^s}^{\Gamma^+} \times J^s$ also induces an
isomorphism of $\mathcal{I}_s^-/\mathcal{I}_s$ onto $\mathcal{K}
(\ell^2 ([0,s]))$.
\end{proof}

\begin{proof}[Proof of \textnormal{(d)}.]
We trivially have $\overline{\bigcup_{r>s} \ker q_r^-} \subset
\ker q_s$. Since $\overline{\bigcup_{r>s} \ker q_r^-}$ is an ideal
in $B_{\Gamma^+} \times_{\tau} \Gamma^+$, there is a
representation $\pi$ of $B_{\Gamma^+} \times_{\tau} \Gamma^+$ on a
Hilbert space $H$ such that $\ker \pi = \overline{\bigcup_{r>s}
\ker q_r^-}$, and then $\pi (i_{\Gamma^+} (t)) = 0$ for $t \geq r
> s$. By Proposition~\ref{univbtaugamma}, there is a
representation $\pi_s$ of $B_{[0,s]} \times_{\tau} \Gamma^+$ on
$H$ such that $\pi_s (i_{\Gamma^+}^{[0,s]} (t)) = \pi
(i_{\Gamma^+} (t))$ for $t \in [0,s]$, and then $\pi_s \circ q_s =
\pi$. Thus $\ker q_s \subset \ker \pi = \overline{\bigcup_{r>s}
\ker q_r^-}$. Now intersecting with $\mathcal{I}$ gives
$\mathcal{I}_s=\overline{\bigcup_{r>s} (\ker q_r^-)\cap
\mathcal{I}}=\overline{\bigcup_{r>s} \mathcal{I}_r^-}$, by, for
example, \cite[Lemma~1.3]{alnr}.
\end{proof}

This completes the proof of Theorem~\ref{Bgammatau}.

\section{The crossed product by the forward shift}\label{sect4}

We now show that when $\Gamma^+$ is the additive semigroup
$\mathbb{N}$, we can obtain more detailed information about the
left-hand and top exact sequences in Theorem~\ref{Bgammatau}.

By viewing functions on $\N$ as sequences, we can identify
$B_{\mathbb{N}}$ with the $C^*$-algebra of convergent sequences
$\mathbf{c}$. Under this identification, the action $\tau$ of
$\mathbb{N}$ is generated by the usual shift $\tau_1: (x_0, x_1,
x_2, \cdots)\mapsto (0, x_0, x_1, x_2, \cdots)$.   The
partial-isometric crossed product
$\mathbf{c}\times_{\tau}\mathbb{N}$ is generated by the power
partial isometry $i_{\N}(1)$, and by \cite[Proposition~9.6]{f} is
the universal $C^*$-algebra generated by a power partial isometry.

Since $[0,n)=[0,n-1]$ for this semigroup, the maps $q_n^-$ and
$q_{n-1}$ coincide, and $\mathcal{I}_n^- = \mathcal{I}_{n-1}$ for
$n>0$. Thus Theorem~\ref{Bgammatau} says that $\mathcal{I}
/\mathcal{I}_0 \cong \mathbb{C}$ and $\mathcal{I}_{n-1} /
\mathcal{I}_n \cong \mathcal{K} (\ell^2 ([0,n])) = M_{n+1}
(\mathbb{C})$ for $n > 0$. We will prove that $\mathcal{I}$ is
isomorphic to $\bigoplus_{n \in \mathbb{N}} M_{n+1} (\mathbb{C})$.
To describe the extensions in the top and left-hand sequences, we
let $P_n:=1 - T^{n+1} (T^*)^{n+1}$ be the projection onto
$\newspan\{e_i:0\leq i\leq n\}$, and define
\[
\mathcal{A} = \big\{f:\mathbb{N} \rightarrow \mathcal{K}
(\ell^2(\mathbb{N})) : f(n) \in P_n \mathcal{K} (\ell^2
(\mathbb{N})) P_n \mbox{ and } \epsilon_\infty(f):=\lim_{n
\rightarrow \infty} f(n) \mbox{ exists } \big\};
\]
note that
$\mathcal{A}_0:=\{f\in\mathcal{A}:\epsilon_\infty(f)=0\}$ is
isomorphic to $\bigoplus_{n \in \mathbb{N}} M_{n+1} (\mathbb{C})$.
Our refinement of Theorem~\ref{Bgammatau} is:

\begin{thm}\label{diagramcNtau}
There are isomorphisms $\pi:\ker \phi_{T^*}\to \mathcal{A}$ and
$\pi^*:\ker \phi_T\to \mathcal{A}$ and an automorphism $\alpha$ of
$\mathcal{A}_0$ such that the following diagram commutes and has
all rows and columns exact:
\begin{equation}\label{comexseqN}
\begin{diagram}\dgARROWLENGTH=0.8\dgARROWLENGTH
\node{} \node{0} \arrow{s} \node{0} \arrow{s} \node{0} \arrow{s}
\\
\node{0} \arrow{e} \node{\mathcal{A}_0} \arrow{e,t}{\id}
\arrow{s,l}{\alpha} \node{\mathcal{A}}
\arrow{e,t}{\epsilon_{\infty}} \arrow{s,l}{\pi^{-1}}
\node{\mathcal{K}(\ell^2(\mathbb{N}))} \arrow{e} \arrow{s}
\node{0}
\\
\node{0} \arrow{e} \node{\mathcal{A}} \arrow{e,t}{(\pi^*)^{-1}}
\arrow{s,l}{\epsilon_{\infty}}
\node{\mathbf{c}\times_{\tau}\mathbb{N}} \arrow{e,t}{\phi_T}
\arrow{s,l}{\phi_{T^*}} \node{\mathcal{T} (\mathbb{Z})} \arrow{e}
\arrow{s,l}{\psi_{T^*}} \node{0}
\\
\node{0} \arrow{e} \node{\mathcal{K}(\ell^2(\mathbb{N}))}
\arrow{e} \arrow{s} \node{\mathcal{T} (\mathbb{Z})}
\arrow{e,t}{\psi_T} \arrow{s} \node{C(\mathbb{T})} \arrow{e}
\arrow{s} \node{0}
\\
\node{} \node{0} \node{0} \node{0}
\end{diagram}
\end{equation}
\end{thm}

Applying the universal property of $\mathbf{c} \times_{\tau}
\mathbb{N}$ to the power partial isometries $P_n T P_n$ and $P_n
T^* P_n$ gives representations $\pi_n$ and $\pi_n^*$ of
$\mathbf{c} \times_{\tau} \mathbb{N}$ on $\ell^2 (\mathbb{N})$
such that $\pi_n (i_{\mathbb{N}} (1)) = P_n T P_n$ and $\pi_n^*
(i_{\mathbb{N}} (1)) = P_n T^* P_n$. We will prove that $\lim_{n
\rightarrow \infty} \pi_n(a)$ exists for all $a\in \ker
\phi_{T^*}$, so that $\pi(a):=\{\pi_n(a)\}$ belongs to
$\mathcal{A}$, and similarly for $\pi^*(b):=\{\pi^*_n(b)\}$ when
$b\in \ker \phi_T$.  To do this, we need to identify  spanning
families for $\ker \phi_{T^*}$ and $\ker \phi_T$.

\begin{lemma}\label{splitexactisom}
For $i,j,m \in \mathbb{N}$, let
\begin{align*}
f_{i,j}^m &:= i_{\mathbb{N}} (i) i_{\mathbb{N}} (m)^*
i_{\mathbb{N}} (m) (1 - i_{\mathbb{N}} (1)
i_{\mathbb{N}} (1)^*) i_{\mathbb{N}} (j)^*,\text{ and }\\
g_{i,j}^m &:= i_{\mathbb{N}} (i)^* i_{\mathbb{N}} (m)
i_{\mathbb{N}} (m)^* (1 - i_{\mathbb{N}} (1)^* i_{\mathbb{N}} (1))
i_{\mathbb{N}} (j).
\end{align*}
Then $\ker \phi_{T^*} = \overline{\newspan} \{f_{i,j}^m: i,j,m \in
\mathbb{N}\}$ and $\ker \phi_T = \overline{\newspan} \{g_{i,j}^m:
i,j,m \in \mathbb{N}\}$.
\end{lemma}

\begin{proof}
First we claim that $\mathcal{E}:= \overline{\newspan}
\{f_{i,j}^m: i,j,m \in \mathbb{N}\}$ is an ideal. To see this, it
suffices to show that $i_{\mathbb{N}} (1) \mathcal{E} \subset
\mathcal{E}$ and $i_{\mathbb{N}} (1)^* \mathcal{E} \subset
\mathcal{E}$. The first is trivial. To show $i_{\mathbb{N}} (1)^*
\mathcal{E} \subset \mathcal{E}$, we first let $i>0$ and compute
using Proposition~\ref{commprojs}:
\begin{align*}
i_{\mathbb{N}} (1)^* i_{\mathbb{N}} (i) &= i_{\mathbb{N}} (1)^*
\big(i_{\mathbb{N}} (1) i_{\mathbb{N}} (i-1)\big) \\
&= i_{\mathbb{N}} (1)^* i_{\mathbb{N}} (1) \big(i_{\mathbb{N}}
(i-1)
i_{\mathbb{N}} (i-1)^* i_{\mathbb{N}} (i-1)\big) \\
&= \big(i_{\mathbb{N}} (i-1) i_{\mathbb{N}} (i-1)^*\big)
\big(i_{\mathbb{N}} (1)^* i_{\mathbb{N}} (1)\big)
i_{\mathbb{N}} (i-1) \\
&= i_{\mathbb{N}} (i-1) i_{\mathbb{N}} (i)^* i_{\mathbb{N}} (i).
\end{align*}
Thus for $i,j,m \in \mathbb{N}$, we have
\begin{align*}
i_{\mathbb{N}} (1)^* f_{i,j}^m &= \left\{ \begin{array}{ll}
i_{\mathbb{N}} (1)^* i_{\mathbb{N}} (i)  i_{\mathbb{N}} (m)^*
i_{\mathbb{N}} (m) \big(1 - i_{\mathbb{N}} (1) i_{\mathbb{N}}
(1)^*\big)
i_{\mathbb{N}} (j)^* & \mbox{if $i > 0$ } \\
i_{\mathbb{N}} (1)^* \big(1 - i_{\mathbb{N}} (1) i_{\mathbb{N}}
(1)^*\big) i_{\mathbb{N}} (m)^* i_{\mathbb{N}} (m) i_{\mathbb{N}}
(j)^* & \mbox{if $ i = 0$}
\end{array}
\right. \\
&= \left\{ \begin{array}{ll} i_{\mathbb{N}} (i - 1)
i_{\mathbb{N}}(i)^* i_{\mathbb{N}} (i) i_{\mathbb{N}} (m)^*
i_{\mathbb{N}} (m) \big(1 - i_{\mathbb{N}} (1) i_{\mathbb{N}}
(1)^*\big) i_{\mathbb{N}} (j)^* & \mbox{if $i > 0$ } \\ 0 &
\mbox{if $i = 0$}
\end{array}
\right. \\
&= \left\{ \begin{array}{ll} i_{\mathbb{N}} (i - 1) i_{\mathbb{N}}
(i \vee m)^* i_{\mathbb{N}} (i \vee m) \big(1 - i_{\mathbb{N}} (1)
i_{\mathbb{N}} (1)^*\big) i_{\mathbb{N}} (j)^* &
\mbox{if $i > 0$ } \\
0 & \mbox{if $i = 0$,}
\end{array}
\right.
\end{align*}
which belongs to $\mathcal{E}$. This proves the claim.

Since $T^*(T^*)^*=1$, each $f_{i,j}^m$ belongs to $\ker
\phi_{T^*}$, and $\mathcal{E} \subset \ker \phi_{T^*}$. Suppose
$\varphi$ is a non-degenerate representation of $\mathbf{c}
\times_{\tau} \mathbb{N}$ on a Hilbert space $H$ with $\ker
\varphi = \mathcal{E}$. Then
\[
1 - \varphi (i_{\mathbb{N}} (1)) \varphi (i_{\mathbb{N}} (1))^* =
\varphi(1 - i_{\mathbb{N}} (1) i_{\mathbb{N}} (1)^*) = 0,
\]
so $\varphi (i_{\mathbb{N}} (1))$ is a coisometry, and $\varphi
(i_{\mathbb{N}} (1))^*$ is an isometry. By Coburn's Theorem, there
is a representation $\psi$ of $\mathcal{T} (\mathbb{Z})$ on $H$
such that $\psi (T) = \varphi (i_{\mathbb{N}} (1))^*$. Then since
$\psi \circ \phi_{T^*}(i_\N(1))= \varphi(i_{\mathbb{N}} (1))$, we
have $\psi\circ\phi_{T^*}=\varphi$, and $\ker \phi_{T^*} \subset
\ker \varphi = \mathcal{E}$.

Similar arguments give the description of $\ker \phi_{T}$.
\end{proof}

\begin{cor}\label{limpipi*}
We have $\pi_n (a)\to \phi_T (a)$ for every $a \in \ker
\phi_{T^*}$, and $\pi^*_n (b) \to\phi_{T^*} (b)$ for every $b \in
\ker \phi_T$.
\end{cor}

\begin{proof}
Since $P_nT^i=0$ unless $i\leq n$, we have
\begin{align}\label{pimform}
\pi_n (f_{i,j}^m) = \pi_n^* (g_{i,j}^m) &= \left\{
\begin{array}{ll}
T^i (1 - TT^*) (T^*)^j & \mbox {if $i,j,m \leq n$} \\
0 & \mbox{otherwise}.
\end{array}
\right.
\end{align}
Since all the homomorphisms have norm $1$, and
$\pi_n(f^m_{i,j})=\phi_T(f^m_{i,j})$ for $n\geq m$, an
$\epsilon/3$-argument shows that $\pi_n (a)\to \phi_T (a)$ for all
$a\in \ker\phi_{T^*}$. Similar arguments give the second
assertion.
\end{proof}

Before proving that $\pi$ and $\pi^*$ map onto $\mathcal{A}$, we
show that they restrict to isomorphisms of
$\mathcal{I}:=(\ker\phi_T)\cap(\ker\phi_{T^*})$ onto
$\mathcal{A}_0$. For this we need a spanning family for
$\mathcal{I}$.

\begin{lemma}\label{Rform}
For $0 \leq i,j \leq m$ we have
\begin{align}\label{fgeq}
f_{i,j}^m - f_{i,j}^{m+1} = g_{m-i,m-j}^m - g_{m-i,m-j}^{m+1},
\end{align}
and the elements \textnormal{(\ref{fgeq})} span $\mathcal{I}$.
\end{lemma}

\begin{proof}
For $i\leq m$, we compute using Proposition~\ref{commprojs}:
\begin{align*}
i_{\mathbb{N}} (i) i_{\mathbb{N}} (m)^* &= i_{\mathbb{N}} (i)
\big(i_{\mathbb{N}} (i)^*
i_{\mathbb{N}} (m-i)^*\big) \\
&= i_{\mathbb{N}} (i) i_{\mathbb{N}} (i)^* \big(i_{\mathbb{N}}
(m-i)^* i_{\mathbb{N}} (m-i)
i_{\mathbb{N}} (m-i)^*\big) \\
&= \big(i_{\mathbb{N}} (m-i)^* i_{\mathbb{N}} (m-i)\big)
\big(i_{\mathbb{N}} (i) i_{\mathbb{N}} (i)^*\big)
i_{\mathbb{N}} (m-i)^* \\
&= i_{\mathbb{N}} (m-i)^* i_{\mathbb{N}} (m) i_{\mathbb{N}} (m)^*.
\end{align*}
Thus for $0 \leq i,j \leq m$, we have
\begin{align*}
f_{i,j}^m - f_{i,j}^{m+1} &= i_{\mathbb{N}} (i)
\big(i_{\mathbb{N}} (m)^* i_{\mathbb{N}} (m) - i_{\mathbb{N}}
(m+1)^* i_{\mathbb{N}} (m+1)\big) \big(1 - i_{\mathbb{N}} (1)
i_{\mathbb{N}} (1)^*\big)
i_{\mathbb{N}} (j)^* \\
&= i_{\mathbb{N}} (i) i_{\mathbb{N}} (m)^* \big(1 - i_{\mathbb{N}}
(1)^* i_{\mathbb{N}} (1)\big) i_{\mathbb{N}} (m) \big(1 -
i_{\mathbb{N}} (1) i_{\mathbb{N}} (1)^*\big)
i_{\mathbb{N}} (j)^* \\
&= i_{\mathbb{N}} (m-i)^* i_{\mathbb{N}} (m) i_{\mathbb{N}} (m)^*
\big(1- i_{\mathbb{N}} (1)^* i_{\mathbb{N}} (1)\big) \cdot\\
&\mbox{\hskip1.5in}\cdot i_{\mathbb{N}} (m-j)\big(i_{\mathbb{N}}
(j) i_{\mathbb{N}} (j)^* -
i_{\mathbb{N}} (j+1) i_{\mathbb{N}} (j+1)^*\big) \\
&= i_{\mathbb{N}} (m-i)^* i_{\mathbb{N}} (m) i_{\mathbb{N}} (m)^*
\big(1 - i_{\mathbb{N}} (1)^*
i_{\mathbb{N}} (1)\big)\cdot\\
&\mbox{\hskip1in}\cdot\big(i_{\mathbb{N}} (m) i_{\mathbb{N}} (m)^*
- i_{\mathbb{N}} (m+1)
i_{\mathbb{N}} (m+1)^*\big)i_{\mathbb{N}} (m-j)\\
&= i_{\mathbb{N}} (m-i)^* \big(i_{\mathbb{N}} (m) i_{\mathbb{N}}
(m)^* - i_{\mathbb{N}} (m+1) i_{\mathbb{N}} (m+1)^*\big)
\cdot\\
&\mbox{\hskip2.5in}\cdot\big(1 - i_{\mathbb{N}} (1)^*
i_{\mathbb{N}} (1)\big)
i_{\mathbb{N}} (m-j) \\
&= g_{m-i,m-j}^m - g_{m-i,m-j}^{m+1}.
\end{align*}

Let $\mathcal{E}_0 = \overline{\newspan} \{f_{i,j}^m -
f_{i,j}^{m+1}: m \in \mathbb{N}, 0 \leq i,j \leq m\}$. Equation
(\ref{fgeq}) implies that $\mathcal{E}_0 \subset \mathcal{I}$.
Since $\ker \phi_T$ and $\ker \phi_{T^*}$ are ideals, $\ker \phi_T
\cap \ker \phi_{T^*} = \ker \phi_T \ker \phi_{T^*}$. A routine
calculation using Proposition~\ref{commprojs} shows that for
$i,j,m,p,r,n \in \mathbb{N}$
\begin{align*}
g_{i,j}^m f_{p,r}^n &= \left\{ \begin{array}{ll}
f_{j+p-i,r}^{j+p} - f_{j+p-i,r}^{j+p+1} & \mbox{if $r,i,m,n\leq j+p$} \\
0 & \mbox{otherwise}
\end{array}
\right.
\end{align*}
which is in $\mathcal{E}_0$. Since $g_{i,j}^m$ and $f_{p,r}^n$
span $\ker \phi_T$ and $\ker \phi_{T^*}$, it follows that $\ker
\phi_T \ker \phi_{T^*}$ is contained in $\mathcal{E}_0$. Thus
$\mathcal{I} = \ker \phi_T \ker \phi_{T^*} = \mathcal{E}_0$.
\end{proof}

\begin{prop}\label{RisomI}
The homomorphisms $\pi:a\mapsto\{\pi_n(a)\}$ and $\pi^*:b\mapsto
\{\pi^*_n(b)\}$ restrict to isomorphisms of $\mathcal{I}$ onto
$\mathcal{A}_0$.
\end{prop}

\begin{proof}
The ideal $\mathcal{A}_0$ is spanned by the functions
$\{e_{ij}^m:0 \leq i,j \leq m\}$ given by
\begin{equation}\label{pinon2}
e_{ij}^m (n) = \left\{ \begin{array}{ll}
T^i (1 - TT^*) (T^*)^j & \mbox{if $m = n$} \\
0 & \mbox{otherwise;}
\end{array}
\right.
\end{equation}
indeed, for each fixed $m$ they span
$P_m\mathcal{K}(\ell^2(\mathbb{N}))P_m$. Equation (\ref{pimform})
implies that
\begin{align}\label{pinonI}
\pi_n (f_{i,j}^m - f_{i,j}^{m+1})  &= \left\{ \begin{array}{ll}
T^i (1 - TT^*) (T^*)^j & \mbox{if $m = n$ and $i,j \leq m$} \\
0 & \mbox{otherwise}
\end{array}
\right.
\end{align}
This proves that $\pi (\mathcal{I})=\mathcal{A}_0$, and similarly
$\pi^* (\mathcal{I})=\mathcal{A}_0$.

The relations (\ref{pinonI}) and (\ref{pinon2}) also show how to
construct an inverse for $\pi$: since
\[
\big\{\{f_{i,j}^m - f_{i,j}^{m+1}:i,j\leq m\}:m\in\mathbb{N}\big\}
\]
consists of mutually orthogonal families of matrix units, there is
a homomorphism of $\bigoplus
M_{m+1}(\mathbb{C})\cong\mathcal{A}_0$ onto $\mathcal{I}$ which
takes $e_{ij}^m$ to $f_{i,j}^m-f^{m+1}_{i,j}$. Similar arguments
using $\mathcal{I}=\overline{\newspan}\{g^m_{i,j}-g^{m+1}_{i,j}\}$
give the corresponding property of $\pi^*$.
\end{proof}

\begin{cor}\label{pipi*onto}
Both $\pi:\ker \phi_{T^*}\to \mathcal{A}$ and $\pi^*:\ker
\phi_{T}\to \mathcal{A}$ are surjective.
\end{cor}

\begin{proof}
Since we know that $\pi(\ker \phi_{T^*})\supset \mathcal{A}_0$, it
suffices to show that for each $K\in \mathcal{K}(\ell^2(\N))$,
there exists $g\in \ker\phi_{T^*}$ with
$\epsilon_\infty(\pi(g))=K$. Indeed, because the range of the
homomorphism $\epsilon_\infty\circ\pi$ is closed, it suffices to
do this for $K= T^i (1 - TT^*)(T^*)^j$. But a computation shows
that $\pi(f_{i,j}^n)=T^i (1 - TT^*)(T^*)^j$ for any $n\geq i\vee
j$. We similarly have $\pi^*(g_{i,j}^{i\vee j})=T^i (1 -
TT^*)(T^*)^j$, and the result follows.
\end{proof}

To see that $\pi$ and $\pi^*$ are injective, we need to know that
$\mathcal{I}$ is an essential ideal in $\mathbf{c} \times_{\tau}
\mathbb{N}$. To do this, we need the following example of a
faithful representation of $\mathbf{c} \times_{\tau} \mathbb{N}$.

\begin{example}\label{excNtau}
Let $V$ be the partial isometry on $\ell^2(\mathbb{N} \times
\mathbb{N})$ such that \[ V(\varepsilon_{k,l}) = \left\{
\begin{array}{ll}
\varepsilon_{k+1,l-1} & \mbox{if $l \geq 1$} \\
0 & \mbox{otherwise,}
\end{array}
\right. \] and note that $V$ is a power partial-isometry, so that
we have a non-degenerate representation $\pi_V \times V$ of
$\mathbf{c} \times_{\tau} \mathbb{N}$ on $\ell^2(\mathbb{N}\times
\mathbb{N})$ such that $\pi_V \times V (i_{\mathbb{N}} (1)) = V$.
If $m>0$ and $i<j$, then
\[
\big(1 - (V^*)^m V^m\big) \big(V^i (V^*)^i - V^j
(V^*)^j\big)(\varepsilon_{0,i})=\varepsilon_{0,i}
\]
so Proposition~\ref{faithgammatau} implies that $\pi_V \times V$
is faithful on $\mathbf{c} \times_{\tau} \mathbb{N}$.
\end{example}

\begin{lemma}\label{Iess}
The ideal $\mathcal{I}$ is essential in $\mathbf{c} \times_{\tau}
\mathbb{N}$.
\end{lemma}

\begin{proof}
Let $V$ be the power partial-isometric representation in
Example~\ref{excNtau}. Then $\pi_V \times V$ is a faithful
representation of $\mathbf{c} \times_{\tau} \mathbb{N}$ on $\ell^2
(\mathbb{N} \times \mathbb{N})$. Since $(\pi_V\times
V)(f^{i+j}_{i,i}-f^{i+j+1}_{i,i})\varepsilon_{i,j}=\varepsilon_{i,j}$,
$\pi_V\times V$ is non-degenerate on $\mathcal{I}$, and it follows
that $\mathcal{I}$ is essential.
\end{proof}

\begin{prop}\label{kerisomA}
Both $\pi : \ker \phi_{T^*} \rightarrow \mathcal{A}$ and $\pi^* :
\ker \phi_T \rightarrow \mathcal{A}$ are isomorphisms.
\end{prop}

\begin{proof}
Corollary~\ref{pipi*onto} says they are surjective. To see that
$\pi$ is injective, suppose $a \in \ker \phi_{T^*}$ and $\pi (a) =
0$. Then for every $c \in \mathcal{I}$, we have $\pi (a c) =0$, $a
c = 0$, and $a = 0$ by Lemma~\ref{Iess}. Thus $\pi$ is injective.
A similar argument shows that $\pi^*$ is injective.
\end{proof}

\begin{proof}[Proof of Theorem~\ref{diagramcNtau}]
We have now proved that we can identify the top and left-hand
sequences with
\begin{equation}\label{exseqCJs}
{0} \longrightarrow {\mathcal{A}_0} \longrightarrow {\mathcal{A}}
\mathop{\longrightarrow}\limits^{\epsilon_\infty}
{\mathcal{K}(\ell^2(\N))} \longrightarrow {0.}
\end{equation}
However, since the isomorphisms $\pi|_\mathcal{I}$ and
$\pi^*|_\mathcal{I}$ are not the same, to make the top left-hand
square commute, we have to introduce an automorphism $\alpha$ of
$\mathcal{A}_0$. The required automorphism is defined on the
spanning elements of (\ref{pinon2}) by
$\alpha(e^n_{i,j})=e^n_{n-i,n-j}$. That the diagram commutes then
follows from Lemma~\ref{Rform}.
\end{proof}

\section{The crossed product by the backward shift}

The backward shifts $\sigma_k$ on $\ell^{\infty}(\mathbb{N})$
satisfy
\[ \sigma_k(1_n) = \left\{ \begin{array}{ll}
1_{n-k} & \mbox{if $n\geq 1$} \\
1 & \mbox{otherwise,}
\end{array}
\right. \] and hence give an action $\sigma:\mathbb{N}\to \End
\mathbf{c}$. In this section we prove a structure theorem for the
crossed product $(\mathbf{c}\times_\sigma \mathbb{N},
k_{\mathbf{c}}, k_{\mathbb{N}})$.

Our first task is to determine the universal property of
$\mathbf{c}\times_\sigma \mathbb{N}$, which is quite different
from that of $\mathbf{c}\times_\tau \mathbb{N}$. The first
difference is that the partial isometries $V^n$ in a covariant
partial-isometric representation $(\pi,V)$ of $(\mathbf{c},
\mathbb{N},\sigma)$ are coisometries:
\[V_n V^*_n = V_n \pi (1) V^*_n =\pi (\sigma_n (1)) = \pi (1) = 1 \mbox{ for
every $n\in \mathbb{N}$}.
\]
The second difference is that the partial isometries
$k_{\mathbb{N}}(n)$ no longer generate $\mathbf{c}\times_\sigma
\mathbb{N}$: we cannot recover $k_{\mathbf{c}}(1_n)$ from
$k_{\mathbb{N}}(n)k_{\mathbb{N}}(n)^*$ alone. More precisely:

\begin{prop}\label{repcoisometry}
Let $(\pi, V)$ be a covariant partial-isometric representation of
$(\mathbf{c},\mathbb{N}, \sigma)$ on $H$, and write $V$ for the
generator $V_1$. Define
\begin{equation}\label{defQ}
Q_0 = 1 - V^*V \ \mbox{ and } \ Q_n := \pi (1_n) - V^* \pi
(1_{n-1}) V \ \mbox{for $n>0$}.
\end{equation}
Then $\{Q_n\}$ is a  sequence of projections satisfying
\begin{equation}\label{sequenceQ}
\cdots \leq Q_{n+1} \leq Q_n \leq Q_{n-1} \leq \cdots \leq Q_0,
\end{equation}
and we can recover $\pi$ via
\begin{equation}\label{picoisom}
\pi (1_n) = (V^*)^n V^n + \sum_{k=0}^{n-1} (V^*)^k Q_{n-k} V^k
\mbox{ for $n>0$.}
\end{equation}
Conversely, for any coisometry $V$ on $H$ and any sequence of
projections $Q_n$ satisfying~\textnormal{(\ref{sequenceQ})}, there
is a covariant partial-isometric representation $(\pi_{V,Q}, V)$
of $(\mathbf{c},\mathbb{N}, \sigma)$ on $H$ such that $\pi_{V,Q}$
satisfies \textnormal{(\ref{picoisom})}.
\end{prop}

\begin{proof}
Suppose $(\pi,V)$ is a covariant partial-isometric representation
of $(\mathbf{c},\mathbb{N},\sigma)$, and define $\{Q_n\}$ using
(\ref{defQ}). Then for $n>0$
\begin{equation}\label{Qproj}
Q_n=\pi(1_n)-V^*\pi(\sigma_1(1_n))V=\pi(1_n)-V^*\big(V\pi(1_n)V^*)V=(1-V^*V)\pi(1_n)
\end{equation}
is the product of commuting projections, and hence is a
projection. We have $Q_1=\pi(1_1)(1-V^*V)\leq 1-V^*V=Q_0$, and for
$n>0$, (\ref{Qproj}) gives
\[
Q_n - Q_{n+1}=(1-V^*V)(\pi (1_n)-\pi (1_{n+1}))\geq 0.
\]
This gives (\ref{sequenceQ}). When we plug the formulas for
$Q_{n-k}$ into the right-hand side of (\ref{picoisom}), the sum
telescopes, and we are left with $\pi(1_n)$.

For the converse, let
\[P_0:=1 \ \mbox{ and }\ P_n:=(V^*)^nV^n+\sum_{k=0}^{n-1}(V^*)^kQ_{n-k}V^k
\mbox{ for $n>0$}.\] Since $Q_n \leq 1 - V^*V$, $VQ_n = 0$. Then
for fixed $n$, $\{(V^*)^k Q_{n-k} V^k: k < n\}$ are mutually
orthogonal and orthogonal to $(V^*)^n V^n$, and hence $P_n$ is a
projection. We have $P_0 - P_1= (1 - V^*V) - Q_1 \geq 0$, and for
$n>0$,
\begin{align*}
P_n - P_{n+1} &= (V^*)^n V^n + \sum_{k=0}^{n-1} (V^*)^k Q_{n-k}
V^k - (V^*)^{n+1} V^{n+1} -
\sum_{k=0}^n (V^*)^k Q_{n+1-k} V^k \\
&= (V^*)^n V^n - (V^*)^{n+1} V^{n+1} - (V^*)^n Q_1 V^n +
\sum_{k=0}^{n-1} (V^*)^k (Q_{n-k} - Q_{n-k+1}) V^k \\
&= (V^*)^n (1 - V^*V - Q_1) V^n + \sum_{k=0}^{n-1} (V^*)^k
(Q_{n-k} - Q_{n-k+1}) V^k,
\end{align*}
and hence $P_n \geq P_{n+1}$. By \cite[Proposition~1.3]{lr}, there
is a representation $\pi_{V,Q}$ of $\mathbf{c}$ such that
$\pi_{V,Q} (1_n) = P_n$.

We now prove that $(\pi_{V,Q}, V)$ is covariant. Since each $V^p$
is a coisometry, we have $\pi_{V,Q} (\sigma_p (1)) = \pi_{V,Q} (1)
= 1 = V^p (V^*)^p$. We also have \[V^p \pi_{V,Q} (1_0) = \pi_{V,Q}
(1_0) V^p = \pi_{V,Q} (\sigma_p(1_0)) V^p.\] For $n>0$, we compute
using $V Q_n = 0$ and Proposition~\ref{commprojs}:
\begin{align*}
V^p \pi_{V,Q} (1_n) &= V^p \Big( (V^*)^n V^n + \sum_{k=0}^{n-1}
(V^*)^k Q_{n-k}
V^k \Big) \\
&= \left\{ \begin{array}{ll} V^p (V^*)^p (V^*)^{n-p} V^{n-p} V^p +
\sum_{k=p}^{n-1} (V^*)^{k-p} Q_{n-k} V^k &
\mbox{if $n > p$} \\
V^{p-n}V^n  + \sum_{k=0}^{n-1}V^{p-k}Q_{n-k}V^k & \mbox{if $n \leq
p$}
\end{array}
\right. \\
&= \left\{ \begin{array}{ll} (V^*)^{n-p} V^{n-p} V^p +
\sum_{k=0}^{n-p-1} (V^*)^k Q_{n-p-k} V^{k+p} &
\mbox{if $n > p$} \\
V^p & \mbox{if $n \leq p$}
\end{array}
\right. \\
&= \left\{ \begin{array}{ll} \big( (V^*)^{n-p} V^{n-p} +
\sum_{k=0}^{n-p-1} (V^*)^k Q_{n-p-k} V^k \big) V^p &
\mbox{if $n > p$} \\
V^p & \mbox{if $n \leq p$}
\end{array}
\right. \\
&= \left\{ \begin{array}{ll}
\pi_{V,Q} (1_{n-p}) V^p & \mbox{if $n > p$} \\
\pi_{V,Q} (1) V^p & \mbox{if $n \leq p$}
\end{array}
\right. \\
&= \pi_{V,Q} (\sigma_p (1_n)) V^p.
\end{align*}
It therefore follows from Corollary~\ref{altcov} that
$(\pi_{V,Q},V)$ is covariant.
\end{proof}

We write $q_n$ for the element
$k_{\mathbf{c}}(1_n)-k_{\mathbb{N}}(1)^*
k_{\mathbf{c}}(\sigma_n(1))k_{\mathbb{N}}(1)$ of
$\mathbf{c}\times_\sigma\mathbb{N}$.
Proposition~\ref{repcoisometry} implies that $\{q_n\}$ is a
decreasing sequence of projections, which together with
$k_{\mathbb{N}}(1)$ generates $\mathbf{c}\times_\sigma\mathbb{N}$.
The next example shows that the $q_n$ are distinct.

\begin{example}\label{egsigma}
Let $V: \mathbb{N} \rightarrow B(\ell^2(\mathbb{N} \times
\mathbb{N}))$ be the coisometric representation such that
\[ V^n (\varepsilon_{k,l}) = \left\{ \begin{array}{ll}
\varepsilon_{k,l-n} & \mbox{if $l \geq n$} \\
0 & \mbox{otherwise,}
\end{array}
\right. \] so that $(V^*)^n (\varepsilon_{k,l}) =
\varepsilon_{k,l+n}$, and let $Q_n$ be the projection on
$\overline{\newspan}\{\varepsilon_{k,0}:k\geq n\}$. Then $\cdots <
Q_{n+1} < Q_n < \cdots < Q_1 < 1 - V^*V$. In the representation
$\pi_{V,Q}$ of Proposition~\ref{repcoisometry}, $\pi_{V,Q}(1_n)$
is the projection on
$\overline{\newspan}\{\varepsilon_{k,l}:k+l\geq n\}$.
\end{example}

We can now characterise the faithful representations of
$\mathbf{c}\times_\sigma\mathbb{N}$:

\begin{prop}\label{coisfaithful}
Suppose $(\pi_{V,Q},V)$ is a covariant partial-isometric
representation of $(\mathbf{c},\mathbb{N},\sigma)$. Then the
representation $\pi_{V,Q} \times V$ of
$\mathbf{c}\times_\sigma\mathbb{N}$ is faithful if and only if
$Q_n\not=Q_{n+1}$ for all $n\geq 0$.
\end{prop}

\begin{proof}

By \cite[Proposition~1.3]{lr}, there is a representation $\pi_Q$
of $\mathbf{c}$ such  that $\pi_Q(1_n) = Q_n$ for $n \in
\mathbb{N}$. For $h \in \range (1 - V^*V) = (V^*H)^{\perp}$, we
have $\pi_{V,Q} (1) h = \pi_Q (1) h$ and
\[
\pi_{V,Q} (1_n) h = \Big( (V^*)^n V^n + \sum_{k=0}^{n-1} (V^*)^k
Q_{n-k} V^k \Big) (1 - V^*V) h = Q_n h = \pi_Q (1_n) h;
\]
thus $\pi_Q = \pi_{V,Q}|_{(V^*H)^{\perp}}$. By
Lemma~\ref{faithful}, $\pi_{V,Q} \times V$ is faithful if and only
if $\pi_{V,Q}$ is faithful on $(V^*H)^{\perp}$, and by
\cite[Proposition~1.3]{lr}, $\pi_{V,Q} = \pi_Q$ is faithful on
$(V^*H)^{\perp}$ if and only if $Q_n \neq Q_{n+1}$ for all $n\geq
0$.
\end{proof}

We are now ready to describe $\mathbf{c}\times_\sigma\mathbb{N}$.
Recall that $T=T_1$ is the unilateral shift on $\ell^2
(\mathbb{N})$, and denote by $F$ the constant function $F:n\mapsto
T^*$. Then $F$ is a coisometry in the $C^*$-algebra $C_b
(\mathbb{N}, B (\ell^2 (\mathbb{N})))$. For $m \in \mathbb{N}$,
define $Q_m \in C_b (\mathbb{N}, B (\ell^2 (\mathbb{N})))$ by
\[
Q_m (n) = \left\{
\begin{array}{ll}
1 - TT^* & \mbox{if $n \geq m$} \\
0 & \mbox{if $n < m$.}
\end{array}
\right.\] Then $\{Q_m\}$ is a decreasing sequence of projections
with $Q_0 = 1 - F^*F$. By Proposition~\ref{repcoisometry}, there
is a homomorphism $\pi_{F, Q} \times F:\mathbf{c} \times_{\sigma}
\mathbb{N}\to C_b (\mathbb{N}, B (\ell^2 (\mathbb{N})))$ such that
$\pi_{F, Q} \times F (k_{\mathbb{N}} (1)) = F$, $\pi_{F, Q} \times
F (q_m) = Q_m$, and
\begin{align*}
\pi_{F, Q} \times F (k_{\mathbf{c}} (1_m)) (n) &= \pi_{F, Q}
\times F \Big(k_{\mathbb{N}} (m)^* k_{\mathbb{N}} (m) +
\sum_{k=0}^{m-1}
k_{\mathbb{N}} (k)^* q_{m-k} k_{\mathbb{N}} (k) \Big) (n) \\
&= \left\{
\begin{array}{ll}
T^m (T^*)^m + \sum_{k=m-n}^{m-1} T^k (1 - TT^*) (T^*)^k & \mbox{if
$n
\leq m$} \\
T^m (T^*)^m + \sum_{k=0}^{m-1} T^k (1 - TT^*) (T^*)^k & \mbox{if
$n > m$.}
\end{array}
\right. \\
&= \left\{
\begin{array}{ll}
T^{m-n} (T^*)^{m-n} & \mbox{if $n \leq m$} \\
1 & \mbox{if $n > m$.}
\end{array}
\right.
\end{align*}

\begin{thm}\label{csigmaN}
The homomorphism $\pi_{F, Q} \times F$ is an isomorphism of
$\mathbf{c} \times_{\sigma} \mathbb{N}$ onto
\[
\mathcal{B} := \{f \in C(\mathbb{N} \cup \{\infty\}, \mathcal{T}
(\mathbb{Z})): \psi_T (f(n)) \mbox{ is constant }\}.
\]
\end{thm}

\begin{proof}
Since $Q_m \neq Q_{m+1}$ for every $m\geq 0$,
Proposition~\ref{coisfaithful} implies that $\pi_{F, Q}\times F$
is faithful. Recall from Proposition~\ref{crossedproduct} that
\[
\mathbf{c} \times_{\sigma} \mathbb{N} =
\overline{\newspan}\{k_{\mathbb{N}} (i)^* k_{\mathbf{c}} (1_m)
k_{\mathbb{N}} (j): i,j,m \in \mathbb{N}\}.
\]
For $i,j,m \in \mathbb{N}$, we have
\begin{equation}\label{pin}
\pi_{F, Q} \times F \big(k_{\mathbb{N}} (i)^* k_{\mathbf{c}} (1_m)
k_{\mathbb{N}} (j) \big) (n) = \left\{
\begin{array}{ll}
T^{i+(m-n)} (T^*)^{j+(m-n)} & \mbox{if $n \leq m$} \\
T^i (T^*)^j & \mbox{if $n > m$.}
\end{array}
\right.
\end{equation}
Thus $\lim_{n\rightarrow\infty} \pi_{F, Q} \times F (a) (n)$
exists for every $a \in \newspan \{k_{\mathbb{N}}
(i)^*k_{\mathbf{c}} (1_m) k_{\mathbb{N}} (j)\}$, and by an
$\epsilon/3$ argument we can extend this to $a\in \mathbf{c}
\times_{\sigma} \mathbb{N}$. Equation (\ref{pin}) also implies
that
\[
\psi_T \big(\pi_{F, Q} \times F (k_{\mathbb{N}} (i)^*
k_{\mathbf{c}} (1_m) k_{\mathbb{N}} (j)) (n)\big) = \epsilon_{i-j}
\]
for every $n \in \mathbb{N}$, and hence $\pi_{F, Q} \times F
(\mathbf{c} \times_{\sigma} \mathbb{N}) \subset \mathcal{B}$.

Let $f \in \mathcal{B}$, and define $g \in \mathcal{B}$ by $g (n)
= f (\infty)$ for every $n \in \mathbb{N}$. Then
\[
\psi_T \big((f-g) (n) \big) = \psi_T (f (n)) - \psi_T (g (n)) =
\psi_T (f (n)) - \psi_T (f (\infty)) = 0,
\]
so $(f-g) (n) \in \ker \psi_T = \mathcal{K} (\ell^2 (\mathbb{N}))$
for all $n \in \mathbb{N}$; since $\lim_{n \rightarrow \infty}
(f-g) (n) = 0$, $f - g$ belongs to $C_0 (\mathbb{N}, \mathcal{K}
(\ell^2 (\mathbb{N})))$. But for $i,j,m \in \mathbb{N}$,
\[
\pi_{F, Q} \times F \big(k_{\mathbb{N}} (i)^* k_{\mathbf{c}} (1_m
- 1_{m+1}) k_{\mathbb{N}} (j) - k_{\mathbb{N}} (i+1)^*
k_{\mathbf{c}} (1_{m-1} - 1_m) k_{\mathbb{N}} (j+1) \big)
\]
is the matrix unit $e_{ij}^m$ of (\ref{pinon2}), so $C_0
(\mathbb{N}, \mathcal{K} (\ell^2 (\mathbb{N}))) =
\overline{\newspan} \{e_{ij}^m: i,j,m \in \mathbb{N}\}$ is
contained in $\pi_{F, Q} \times F (\mathbf{c} \times_{\sigma}
\mathbb{N})$.  The function $g$ is constant, so it belongs to
\[
C^*(F)=\pi_{F, Q} \times F (C^* (k_{\mathbb{N}} (1))) \subset
\pi_{F, Q} \times F (\mathbf{c} \times_{\sigma} \mathbb{N}),
\]
and hence so does $f=(f-g)+g$.
\end{proof}

\begin {cor}\label{exseqcsigma}
There is an exact sequence
\begin{equation*}
0\longrightarrow{C (\mathbb{N} \cup \{\infty \}, \mathcal{K}
(\ell^2 (\mathbb{N})))} \longrightarrow{\mathbf{c} \times_{\sigma}
\mathbb{N}}\mathop{\longrightarrow} {C
(\mathbb{T})}\longrightarrow 0.
\end{equation*}
\end{cor}

\end{document}